\documentclass[10pt]{article}
\usepackage[T1]{fontenc} 
\usepackage[utf8]{inputenc}

\usepackage{lmodern}

\usepackage{amssymb} 
\usepackage{latexsym}
\usepackage{epsfig}
\usepackage{lscape,graphicx,color,psfrag}
\usepackage{calrsfs}
\usepackage[cmex10]{amsmath}
\usepackage{amsthm}

\theoremstyle{plain}

\newtheorem{proposition}{Proposition}

\newtheorem{conjecture}{Conjecture}
\theoremstyle{definition}

\theoremstyle{remark}
\newtheorem{remark}{Remark}

\def\limsup{\mathop{\textup{lim\,sup}}}

\def\minim{\mathop{\textup{minimize}}}
\def\st{\mathop{\textup{subject to}}}

\newcommand{\Ex}{\ensuremath{\mathsf{E}}}

\usepackage{appendix}

\usepackage[top=1in, bottom=1.25in, left=0.9in, right=0.9in]{geometry}

\usepackage[sort&compress]{natbib} 
\setcitestyle{authoryear,open={(},close={)}}

\usepackage[colorlinks=true,breaklinks=true,bookmarks=true,urlcolor=blue,
     citecolor=blue,linkcolor=blue,bookmarksopen=false,draft=false]{hyperref}

\def\EMAIL#1{\href{mailto:#1}{#1}}
\def\URL#1{\href{#1}{#1}}         

\begin{document}

\title{Resource allocation and routing in parallel multi-server queues with abandonments for cloud profit maximization}

\author{Jos\'e Ni\~no-Mora \\ Department of Statistics \\
    Carlos III University of Madrid \\
     28903 Getafe (Madrid), Spain \\  \EMAIL{jose.nino@uc3m.es}, \URL{http://orcid.org/0000-0002-2172-3983} }
     
\date{Published in \textit{Computers \& Operations Research}, vol.\ 103, 221--236, 2019 \\
DOI: \href{https://doi.org/10.1016/j.cor.2018.11.012}{10.1016/j.cor.2018.11.012}
}

\maketitle

\begin{abstract}
This paper considers a Markov decision model for profit maximization of a cloud computing service provider catering to customers submitting jobs with firm real-time random deadlines. 
Customers are charged on a per-job basis, receiving a full refund if deadlines are missed.
The service provider leases computing resources from an infrastructure provider 
in a two-tier scheme: long-term leasing of basic infrastructure, consisting of heterogeneous parallel service nodes, each modeled as a multi-server queue, and short-term leasing of external servers. 
Given the intractability of computing an optimal dynamic resource allocation and job routing policy,  maximizing the long-run average profit rate,
the paper addresses the design, implementation and testing of low-complexity heuristics.  
The  policies considered are  
a static policy given by an optimal Bernoulli splitting, and three dynamic index policies based on different index definitions: individually optimal (IO), policy improvement (PI) and 
restless bandit (RB) indices. 
The paper shows how to implement efficiently  each such policy, and presents a comprehensive empirical comparison, drawing qualitative insights on their strengths and weaknesses, and benchmarking their performance in an extensive study.
\end{abstract}



\section{Introduction} 
\label{s:intro} 
\subsection{Motivation}
\label{s:rarrccp}
To keep pace with steeply growing global demand, 
e-service providers
increasingly rely on platforms integrating heterogeneous computing resources. 
Such is the case with \emph{cloud computing}, which has emerged over the last decade aiming to realize the vision of computing as the 5th utility.
See, e.g., \citet{buyya09}, \citet{li13}, \citet{caoLiStoj14}, and \citet{meietal15}.

A cloud computing environment involves three main stakeholders: 
customers, service providers, and infrastructure providers.
In exchange for fees, 
customers expect to receive a certain Quality of Service (QoS) level. 
Penalties to the service provider for degraded QoS are specified in a
Service Level Agreement (SLA), which may include clauses for  refund of service fees.

Service providers need to decide how to provision computing resources
to maximize profit.  
Instead of owning such resources,  it is often more economical to lease them from infrastructure providers.
In a \emph{static resource allocation} scheme, a fixed 
\emph{basic infrastructure} is leased on a long-term basis, consisting of heterogeneous parallel multi-server nodes.
In a \emph{dynamic resource allocation} scheme (see \citet{meietal15}) such basic infrastructure is complemented with extra external resources that are occasionally leased on a short-term basis when deemed convenient, e.g., due to temporary overload.

To fully specify how incoming requests, thereafter referred to as \emph{jobs}, are handled in the latter setting, 
the service provider needs to adopt a \emph{joint resource allocation and job routing policy}, which
prescribes for each incoming job whether to route it to a service node in the basic infrastructure, or to outsource it instead to an external server.
Since the average profit depends on the policy adopted, this motivates the research goal of designing policies that are simple to implement and yet nearly optimal for maxmizing profit.

This paper addresses such an issue as it applies to emerging cloud platforms providing time-critical services,    catering
 to impatient customers whose 
 QoS requirements take the form of \emph{random firm deadlines}, which are unknown by the service provider until they expire. 
In such platforms, jobs immediately abandon upon missing their deadlines, as they lose all value to the customers.
 See, e.g., \citet{pdt14} and \citet{chiangetal16}.
 
We consider the two standard types of firm deadlines: 
 \emph{deadlines to the beginning of service} (DBS) and \emph{deadlines to the end of service} (DES). 
A job's deadline specifies, under DBS,  that its service should begin before a certain
time, and, under DES, that it should end before a certain
time. 
Thus, an incoming job with a \emph{relative deadline} (time from arrival to deadline expiration) of 5 min.\ should start service within 5 min.\ under DBS, and complete service within 5 min.\ under DES. Otherwise, it abandons.

Examples are found, e.g.,  in \emph{distributed real-time database} applications (see \citet{lasotaetal17}), in particular in those where content is replicated across multiple servers and transactions must be started or completed before current conditions change significantly. Think 
of online low-latency high-frequency stock trading platforms (see  \citet{hasbrouckSaar13}), aiming to complete transactions before current stock prices move beyond given limits.

Another example is 
 \emph{online retailing} platforms. 
Since online shoppers are willing to wait only a limited time for a page to load before navigating away (see   \cite{nah04}, and \citet{priyaetal17} on shopping cart abandonment), it is important for such platforms to minimize the fraction of lost customers.

An emerging firm real-time application is numerical weather prediction (see \citet{siutaetal16}), where 
the customer is a weather forecaster submitting computationally demanding models to run on a cloud platform, with current weather conditions as input. If the latter change before the model run completes, the job loses all value and is dropped. 

 In firm real-time platforms, jobs that end up abandoning are harmful not only in that  
they contribute no profit (if the service fee is refunded), but because their sojourn in the system causes later jobs to also 
abandon, further lowering profit.
A standard approach is to 
 incorporate \emph{admission control} so that jobs can be rejected on arrival, as in   \citet{wuetal12} and \citet{chiangetal16}.
However, upfront job rejection has undesirable effects, such as a sure 
 loss of both revenue and customer goodwill. 
An alternative is to use dynamic resource allocation as explored herein.

\subsection{Model formulation}
\label{s:md}
\begin{figure}[htb!]
    \centering
    \psfrag{a}{$\lambda$}
    \psfrag{b}[cb]{$A(t) = 0$}
    \psfrag{c}[c]{$A(t) = 1$}
    \psfrag{d}[lb]{$A(t) = n$}
    \psfrag{g}[rb]{$L_1(X_1(t))$}
    \psfrag{h}[lb]{$L_n(X_n(t))$}
    \psfrag{e}[cb]{$\mu_1$}
    \psfrag{f}[cb]{$\mu_n$}
    \psfrag{i}[cb]{$\mu_0$}
\includegraphics[width=0.75\textwidth]{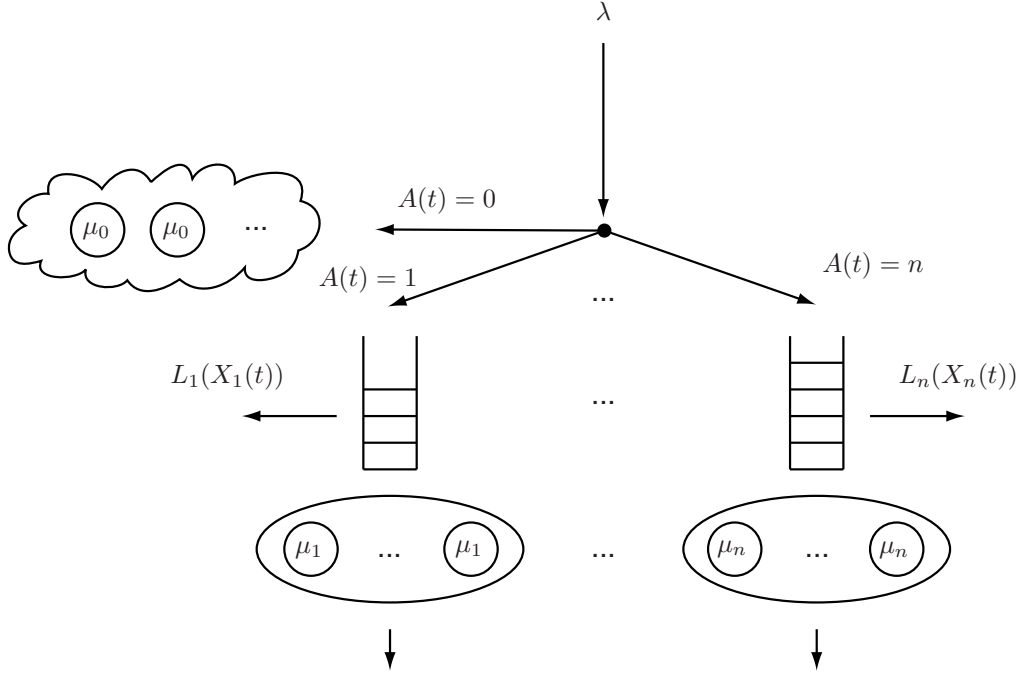}
    \caption{Dynamic resource allocation and job routing model with abandonments.}
    \label{fig:acrmodel}
\end{figure}

This paper considers a
\emph{Markov decision process} (MDP) model (see \citet{put94}) of a cloud platform as outlined above. See  Fig.\ \ref{fig:acrmodel}. 
Basic long-term leased infrastructure is modeled as a collection of $n$ parallel multi-server nodes, with
service node $k = 1, \ldots, n$ 
having its own queue with unlimited buffer and a finite number $m_k$ of identical exponential servers with   rate $\mu_k$.
As for external short-term leased infrastructure, it is modeled as a multi-server node labeled by $k = 0$ with $m_0 \triangleq \infty$ servers with rate $\mu_0$.  
Jobs arrive as
a Poisson stream with 
rate $\lambda$ and independent service times.

The \emph{relative deadlines} (time from arrival to deadline expiration) of jobs are modeled as independent exponential random variables with \emph{abandonment rate} $\theta$,  being also independent of arrival and service times. 
The \emph{total abandonment} (or \emph{loss}) \emph{rate} $L_k(i_k)$ for node $k$ when it lies in \emph{state} (number of jobs present)
$X_k(t) = i_k \in \mathbb{Z}_+ \triangleq \{0, 1, \ldots\}$  is $L_k(i_k) \triangleq (i_k-m_k)^+ \theta$ under DBS, where $x^+ \triangleq \max(x, 0)$, and $L_k(i_k) \triangleq  i_k \theta$ 
under DES.
In \citet{baccelliHebu81}'s notation, node $k$
is modeled  as an M/M/$m_k+$M queue.

Immediately upon a job's arrival at time $t$, 
the system controller sends it irrevocably to a service node
 $A(t) \in \{0, 1, \ldots, n\}$, where such \emph{actions} are prescribed through a
policy $\boldsymbol{\pi}$ from the class $\boldsymbol{\Pi}$ of nonanticipative policies.
Thus, action $A(t)$ 
is a (possibly randomized) function of  the \emph{observed history} 
$\mathcal{H}(t) \triangleq 
\{(\mathbf{X}(s), A(s))\colon 0 \leqslant s < t\} \cup \{\mathbf{X}(t)\}$
 of states and actions, where
the \emph{system state} at time $t$ is  
$\mathbf{X}(t) \triangleq (X_k(t))_{k=1}^n$. Note that the number of jobs $X_0(t)$ at the external node is not considered here part of the system state.

Jobs dispatched to a basic node are scheduled in \emph{first-come first-serve} (FCFS) order, which
 ensures that the sojourn-time distribution of each such job is determined by the node's state found upon arrival. 

The service provider charges customers an upfront
 service fee of $F > 0$ per job,  normalized here to $F = \$1$, which is fully refunded if the job's deadline is missed so it abandons.
 For each job outsourced to the external node, a lump 
 charge $c> 0$ is paid to the infrastructure provider.
 Note that the latter jobs may still miss their deadlines under DES,  with probability 
 $\theta/(\theta + \mu_0)$. 
 Thus, the total expected cost $C$ of allocating a job to the external node, including the possible refund,  is 
 $C \triangleq c$ under DBS and 
 $C \triangleq c +  \theta/(\theta + \mu_0)$ under DES.
 We shall assume that $C < 1$, as otherwise it would clearly be uneconomical to use the external node.
 
In such a setting, this paper addresses the design and implementation of effective joint dynamic resource allocation and job routing policies, aiming to maximize the service provider's long-run average expected profit per unit time.
Since the average rate of service fees collected per unit time equals $\lambda$, such an objective is equivalent to minimizing the long-run average 
expected rate per unit time of deadline-missing refunds and short-term leasing charges for use of the external node.
The latter problem can be formulated as

\begin{equation}
\label{eq:panew}
\minim_{\boldsymbol{\pi} \in \boldsymbol{\Pi}} \, \limsup_{T \to \infty} \frac{1}{T}
\, \Ex_{\mathbf{i}}^{\boldsymbol{\pi}}\Big[\int_0^T  \Big(\sum_{k =1}^n
  L_k(X_k(t)) + \lambda  C 1_{\{A(t) = 0\}}\Big) \, dt\Big],
\end{equation}
where $\Ex_{\mathbf{i}}^{\boldsymbol{\pi}}$ denotes expectation starting from $\mathbf{X}(0) = \mathbf{i} = (i_k)_{k=1}^n$ under policy $\boldsymbol{\pi}$ and ``$\limsup$'' denotes limit superior.
A policy $\boldsymbol{\pi}^* \in \boldsymbol{\Pi}$ is \emph{average-cost optimal} if it  minimizes the objective in (\ref{eq:panew})  for every  possible initial state $\mathbf{i}$.

In the real-time environment of concern here, a practical requirement on policies is that they 
 allow a low-complexity implementation preventing burdensome overheads. 
 Note that, in the computer communications literature, use of dynamic routing policies basing decisions on the current system state is often considered impractical,  due to the communication latency for gathering distributed state information. See 
 \citet{heetal06}.
Yet,  dynamic policies might be practically implementable for  cloud computing applications where  latency is negligible relative to job processing times and deadline durations. Think, e.g., of the numerical weather prediction application mentioned above, where processing of each job can take minutes or hours.

Problem (\ref{eq:panew}) is a denumerable state MDP with unbounded state transition and cost rates.
Under certain conditions, reviewed in \citet[Ch.\ 7]{guoHL09}, such MDPs have 
optimal \emph{stationary deterministic} policies, which select actions based on the current state,  characterized by the problem's  \emph{dynamic programming} (DP) equations. 
Yet, 
computing an optimal policy for the above model by solving numerically its DP equations is generally not possible, as there is an infinite number of equations, one per state.
Even if one considers a finite-state model with finite-buffer queues, 
the DP equations cannot be solved for platforms with more than a few nodes due to the \emph{curse of dimensionality}, as the number of states grows exponentially with the number of nodes. 

In \S\ref{s:obsei} and \S\ref{s:civdpni} we shall refer to 
three base instances of the above model with $n = 3$ basic nodes, both under DBS and DES, 
with parameters as shown in Table \ref{t:bi}, where $\mathbf{m} = (m_1, m_2, m_3)$ and $\boldsymbol{\mu} = (\mu_1, \mu_2, \mu_3)$.
Note that base instance $1$ represents a balanced system in that basic nodes have equal total processing capacities $m_k \mu_k \equiv 20$.
Base instance $2$ represents an imbalanced system with 
total processing capacities of basic nodes 
ordered as $m_1 \mu_1 < m_2 \mu_2 < m_3 \mu_3$, whereas 
base instance $3$ has $m_2 \mu_2 < m_1 \mu_1 < m_3 \mu_3$.
The arrival rate $\lambda$ has been chosen so that 
the \emph{nominal system load} $\rho \triangleq \lambda / \sum_{k=1}^3 m_k \mu_k$ equals $1$ in each base instance. 
Note that the information given on the external node is reduced to its expected usage cost $C$, as there is no need to specify $c$ or $\mu_0$.
The job's abandonment rate is $\theta = 0.3$, so the mean relative deadline is $10/3$  time units.

\begin{table}[!htb]
\centering
\begin{tabular}{|l|cccccc|}
\hline
 & $n$ & $\lambda$ & $\theta$ & $C$ & $\mathbf{m}$ & $\boldsymbol{\mu}$\\ \hline
Base instance 1: & $3$ & $60$ &  $0.3$ & $0.2$ & $(2, 5, 10)$ & $(10, 4, 2)$ \\
Base instance 2: & $3$ & $53.5$ &  $0.3$ & $0.2$ & $(2, 5, 10)$ & $(8, 3.5, 2)$ \\
Base instance 3: & $3$ & $62.5$ &  $0.3$ & $0.2$ & $(2, 5, 10)$ & $(10, 3.5, 2.5)$ \\
\hline
\end{tabular}
\caption{Base instances for \S\ref{s:obsei} and \S\ref{s:civdpni}.}
\label{t:bi}
\end{table}

\subsection{Heuristic policies based on Bernoulli splitting and routing indices}
\label{s:fadhp}
The above discussion motivates the interest of designing heuristic policies that, though suboptimal, can be implemented with low complexity  and perform well.
For such a purpose, this paper deploys 
four approaches.

The first approach produces a \emph{static policy} (blind to the state) 
given by 
a \emph{Bernoulli splitting} (BS) of the arrival stream, where the total arrival rate $\lambda$ is split into 
a rate $\lambda_k$ for each node $k = 0, 1, \ldots, n$, and then each arrival is sent to node $k$ with probability $\lambda_k/\lambda$, independently of other jobs. 
In an \emph{optimal BS} the rates $\lambda_k$ are chosen to maximize average profit \emph{within the class of BS policies} (so they are \emph{optimal} in that limited sense). See  \citet{lee95}.

The other three approaches yield dynamic \emph{index policies} that use the system state in a 
 tractable way, being based on 
 \emph{routing indices} $\varphi_k(i_k)$ attached to each basic node $k = 1, \ldots, n$ as a function of its state $i_k$ (number of jobs present), which are used as a measure of \emph{routing priority}: the lower the index value of a node the higher the priority for routing a job to it.
Thus, upon a job's arrival,
 it is dispatched to a basic node with currently \emph{lowest index value} (breaking ties arbitrarily), provided that the latter does not exceed the external node's expected usage cost $C$. Otherwise, the job is outsourced to the external node.
If the indices  $\varphi_k(i_k)$ can be computed with low complexity, and if 
in the application at hand the latency for 
gathering the system state is negligible, 
  such policies may be suitable for practical implementation.

For such a policy to make intuitive sense,
 the index $\varphi_k(i_k)$
must somehow measure the cost of routing a job to basic node $k$ when this lies 
in state $i_k$. 
This paper considers 
three approaches to gauging such routing costs.

The first approach measures the cost of routing a job to basic node $k$ by considering only the impact on that particular job.  
This leads to taking as $\varphi_k(i_k)$ the expected refund from the job if routed to node $k$, which yields the
\emph{individually optimal} (IO) \emph{index policy}.

The second approach takes one step of the \emph{policy improvement} (PI) algorithm for MDPs starting from an optimal BS, 
yielding the \emph{PI index policy}.
 
The third approach casts the model into the framework of the \emph{multi-armed restless bandit problem} (MARBP), and then deploys \citet{whit88b}'s index policy.
This yields the \emph{restless bandit} (RB) \emph{index policy}.

\subsection{Target properties for routing indices}
\label{s:dpeop}
A basic node's routing index $\varphi(i)$ (where the node label $k$ is dropped from the notation) is a function of the node's state $i$, which somehow incorporates some model parameters.
Based on intuitive grounds, it appears reasonable to consider routing indices that satisfy certain structural properties.
In particular,  the target properties listed in Table \ref{t:indexp} are proposed herein, concerning dependence on the node's state $i$, the system's arrival and abandonment rates $\lambda$ and $\theta$, the external node's (expected) usage cost $C$, and the node's servers' rate $\mu$ and server pool size $m$.
Below we use the acronyms \emph{nondecreasing} (ND) and \emph{nondecreasing} (NI).

\begin{table}[!htb]
\centering
\begin{tabular}{|ll|}
\hline
 & Other things being equal,
 \\
$P_1\colon$ & $\varphi(i)$ is ND in $i$ with 
$\varphi(i) \nearrow 1$ as $i \nearrow \infty$. \\
$P_2\colon$ & $\varphi(i)$ is ND in $\lambda$ \\
$P_3\colon$ & $\varphi(i)$ is ND in $\theta$, with $\varphi(i) \searrow 0$ as $\theta \searrow 0$ and
$\varphi(i) \nearrow 1$ as $\theta \nearrow \infty$ \\
$P_4\colon$ &   $\varphi(i)$ is NI in $C$ for large enough $C$ \\
 $P_5\colon$ & $\varphi(i)$ is NI in  $\mu$, with $\varphi(i) \nearrow 1$ as $\mu \searrow 0$ and
$\varphi(i) \searrow 0$ as $\mu \nearrow \infty$ \\
$P_6\colon$ &  $\varphi(i)$ is NI in $m$ \\
\hline
\end{tabular}
\caption{Target properties for a routing index $\varphi(i)$.}
\label{t:indexp}
\end{table} 

\begin{remark}
\label{re:tpri}
The following is an intuitive rationale for properties in Table \ref{t:indexp}. 
\begin{itemize}
\item[\textup{(i)}] $P_1$: the more congested a basic  node, the lower its routing priority, which can become lower than that of the external node.
\item[\textup{(ii)}] $P_2$: the higher the arrival rate, 
the lower the routing priorities of basic nodes relative to the external node.
\item[\textup{(iii)}] $P_3$: the more impatient the customers,  the lower the routing priorities of basic nodes relative to the external node.
\item[\textup{(iv)}] $P_4$: the routing priority of a basic node does not increase with $C$, for large enough $C$.  
\item[\textup{(v)}] $P_5$: the faster a node's servers, the higher its routing priority. 
\item[\textup{(vi)}] $P_6$: the larger a node's server pool, the higher its routing priority. 
\end{itemize}
\end{remark}

The above raises the question of whether properties $P_1$--$P_6$ in Table \ref{t:indexp} are consistent with structural properties of optimal policies, as intuition would suggest. 
Yet, resolving such an issue is beyond the scope of this paper. 

\subsection{Goals and contributions}
\label{s:gandc}
The main goal of this paper is to compare the approaches mentioned above to the design of policies for the present model, with respect both to the complexity of implementing them and to their empirical  performance. 

The paper aims to assess the effects of changes in model parameters on the deviation from optimality of such  policies to identify their strengths and weaknesses.
For such a purpose, an extensive numerical study is carried out on a wide range of instances. 

Contributions include the following: 
(1) formulation of a new MDP model for profit maximization of a cloud platform for impatient customers with dynamic resource allocation and job routing; 
(2) development of four tractable heuristic policies, along with efficient means for their computation;
 the means of computing an optimal BS is based on an empirically supported conjecture
proposed here on properties of performance metrics for the M/M/$m+$M queue, as well as on new relations reducing the analysis of the M/M/$m+$M queue under DES to the DBS case; 
(3) identification of qualitative insights on the heuristics considered; and  
(4) an extensive comparative numerical study on the performance of the policies across a wide range of instances. 

\subsection{Organization of the paper}
\label{s:oofp}
The remainder of the paper is organized as follows. 
\S \ref{s:rw} reviews related work.
\S \ref{s:obsm} considers the BS policy.
\S \ref{s:ip} discusses the IO, PI and RB policies.
\S \ref{s:cbs} reports  the results of a comparative benchmarking study on the performance of the policies considered.
\S \ref{s:d} presents a final discussion of results.
Two appendices contain ancillary material.
\ref{s:evalell} presents required results for computing the optimal BS on analysis of the M/M/$m+$M queue under DBS, and shows how to reduce the analysis of the DES case to the DBS case.
\ref{a:rbip} outlines how to reformulate the present model into the framework of the MARBP, and reviews computation of the RB policy. 

\section{Some related work}
\label{s:rw}
The study of multi-dimensional MDP models for optimal dynamic resource allocation has attracted substantial research attention. 
Since the numerical solution of their DP equations is hindered by the curse of dimensionality, researchers have sought to identify optimal policies with a simple structure, often of index type. 
See, e.g., \citet{courVa83}, \citet{katehDerm84}, and \citet{katehJoh84}, which are among the first papers to use MDPs for
this kind of problems.
Further optimality results are obtained in 
\citet{katehMel88,katehMel95},  where the latter paper also considers optimal routing to queues.

For more complex models, identifying optimal policies is elusive. Yet, works such as \citet{katehLev86} and \citet{katehDerm89} presented methods
to analyze the models investigated under light traffic and heavy traffic conditions, giving efficient algorithms for finding optimal or asymptotically policies in such regimes.

The design and implementation of resource allocation and/or job routing policies in MDP models of cloud or similar e-service platforms has been the subject of substantial recent research attention.
Markovian multi-server queues have been often used as system models in such settings: for blade servers in \citet{li13}, multi-core computer servers in \citet{caoLiStoj14}, and multi-server cloud platforms in \citet{meietal15}.
The latter work considered an M/M/$m+$D model of a cloud service provider where all jobs have the same relative deadline.
That work assumes a dynamic resource allocation scheme where jobs missing deadlines are  sent to a short-term leased server.

Research on queueing models with 
abandonments originated with the work of \citet[60--67]{palm57b} on the M/M/$m+$M queue under DBS, and has  been extensively developed. 
See \citet{baccelliHebu81}.
Currently, such models are mostly applied to the analysis of call centers. See, e.g.,  \citet{zeltynMan05}.
More recent work applies them to the study of cloud platforms, as in \citet{pdt14} and \citet{chiangetal16}.

Optimal BS policies have been developed for several queueing models, but not for that herein.
\citet{buzenChen74} derived the optimal BS for 
routing jobs to parallel M/G/$1$ queues to minimize mean response time.
\citet{lee95} extended such work to a model where jobs are classified into multiple priority classes.
\citet{heetal06} applied optimal BS to a model with parallel 
M/M/$m$ queues. 
\citet{li13} obtained the optimal BS for routing generic jobs to parallel  M/M/$m$ queues that also cater to dedicated jobs.
\citet{kallmescass95} applied optimal BS for deadline-miss rate minimization in a model with admission control and routing of soft real-time jobs to parallel M/M/$1$ queues.
\citet{heetal06} extended the latter work to parallel M/M/$m$ queues, and 
\citet[\S 4.1]{nmcor12} further incorporated admission control.

As for IO index policies, they are optimal in certain models for routing jobs to parallel symmetric queues 
(see \citet{winston77} and \citet{johri89} for mean response time minimization, and \citet{movagh05} for deadline-miss rate minimization with firm real-time jobs).
Use of such policies in heterogeneous systems has been addressed, e.g., 
 in \citet{chowKohler79} for mean response time minimization with parallel M/M/$1$ queues, and in \citet[\S 3]{nmcor12} for deadline-miss rate minimization of soft real-time jobs with parallel M/M/$m$ queues  and admission control.

Work on PI index policies includes, e.g.,  \citet{krish90} and \citet{sassen97}, which considered minimization of mean job response time with parallel M/M/$m$ queues and M/G/$1$ queues, respectively. 
\citet[\S 4.2]{nmcor12} developed a PI index policy for control of admission and routing to parallel M/M/$m$ queues with soft real-time jobs for deadline-miss rate minimization.
 
Regarding Whittle's RB index policy, its application to 
admission control and routing to parallel queues was introduced in \citet[\S 8.1]{nmmp02}, in a broad model including that herein as a special case. See also \citet{nmnetcoop07}. 
\citet[\S 5]{nmcor12} investigated such a policy in a model for control of admission and routing of soft real-time jobs to parallel M/M/$m$ queues.  

\section{Optimal BS policy}
\label{s:obsm}
{This section develops a heuristic static policy for problem (\ref{eq:panew}) that is optimal within the class of BS policies, adapting the approach in 
\citet{buzenChen74} to the present model.
In a BS, the actions selected upon job arrivals are drawn according to fixed probabilities:
when a job arrives at time $t$, it is dispatched to service node $k = 0, 1, \ldots, n$ (i.e., 
action $A(t) = k$ is selected) 
with probability $p_k = \lambda_k/\lambda$, independently of the actions taken on previous jobs, where the $\lambda_k$'s are rates adding up to $\lambda$ to be determined. 

The input to service node $k$ under such a BS policy is a Poisson process with rate $\lambda_k$, and hence 
the node behaves as an M/M/$m_k+$M queue with \emph{offered load} $r_k(\lambda_k) \triangleq 
\lambda_k/\mu_k$ and \emph{offered load per server} $\rho_k(\lambda_k) \triangleq 
r_k(\lambda_k)/m_k$, which
is stable.
Let $P_{k, \mathrm{ab}}(\lambda_k)$ be the \emph{abandonment probability} for node $k$, i.e., the probability that a random arrival abandons due to missing its deadline. 

An \emph{optimal BS} is a globally optimal solution $\boldsymbol{\lambda}^* = (\lambda_k^*)_{k=0}^n$ to 
the following constrained optimization problem (cf.\ (\ref{eq:panew})), where $\ell_k(\lambda_k) \triangleq \lambda_k P_{k, \mathrm{ab}}(\lambda_k)$ is the \emph{mean abandonment} (or \emph{loss}) \emph{rate} for basic node $k = 1, \ldots, n$:
\begin{equation}
\label{eq:nlp}
\begin{split}
& \minim \, \sum_{k =1}^n \ell_k(\lambda_k) + C \lambda_0 \\
& \st\colon \lambda_0, \lambda_1, \ldots, \lambda_n \geqslant 0 \\
& \sum_{k =0}^n \lambda_k = \lambda.
\end{split}
\end{equation}

The evaluation of functions $\ell_k$ and their derivatives $\ell_k'$, which is required 
to solve numerically problem (\ref{eq:nlp}), is addressed in \ref{s:evalell}.

\subsection{Computing the optimal BS}
\label{s:cobs}
We next address computation of the optimal BS.
Since (\ref{eq:nlp}) is a 
smooth linearly constrained nonlinear optimization problem, 
standard results (see, e.g., \citet{andreasApr11}) yield that, if $\boldsymbol{\lambda}^*$
is a
\emph{local optimum}, there is a Lagrange multiplier $\alpha^*$ for its equality
constraint satisfying the  
 \emph{Karush--Kuhn--Tucker} (KKT) \emph{first-order optimality conditions}:
\begin{equation}
\label{eq:kktc}
\begin{split}
C & \geqslant \alpha^*, \quad \text{ with ``='' if } \lambda_0^* > 0 \\
\ell_k'(\lambda_k^*)  & = \alpha^* \quad \text{ for any basic node $k$ having }  \lambda_k^* > 0 \\
\alpha_k  & \geqslant \alpha^* \quad \text{ for any basic node $k$ having } \lambda_k^* = 0,
\end{split}
\end{equation}
where $\alpha_k \triangleq \ell_k'(0^+) = P_{k, \mathrm{ab}}(0)$ is the probability that an arrival to an empty node $k$ abandons, so
 $\alpha_k = 0$ under DBS and $\alpha_k = \theta/(\theta + \mu_k)$ under DES.

To ensure that such conditions uniquely determine a global optimum for problem (\ref{eq:nlp})
we need the functions $\ell_k$ in the objective to satisfy certain  properties such as, e.g.,  being strictly convex. Yet, 
establishing such properties is a challenging analytical problem beyond the scope of this paper, 
due to the 
complexity of the formulae involved (see \ref{s:evalell}), which, to the author's knowledge,  has not been addressed hitherto in the literature. 

Still, based on observation of particular instances (see, e.g., Fig.\ \ref{fig:ffprime}) the following  conjecture on properties of the mean abandonment rate
 $\ell(\lambda)$ for an M/M/$m+$M queue is proposed herein, where the  label $k$ is dropped from the notation. 
Note that the term ``increasing'' is used below in the strict sense.

\begin{figure}[!ht]
\centering
\includegraphics[height=1.7in]{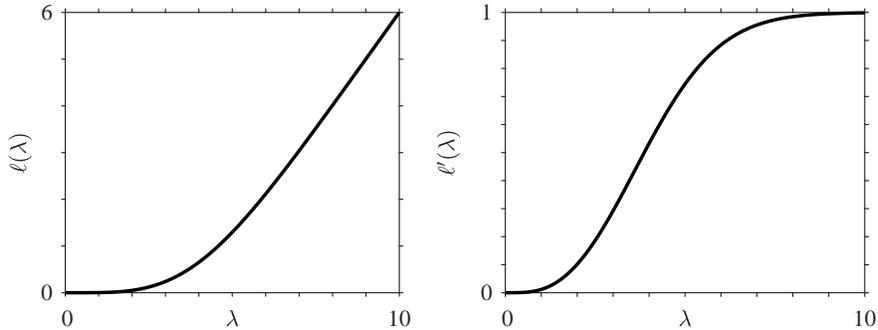}
\caption{$\ell(\lambda)$ and $\ell'(\lambda)$ for an M/M/$m+$M queue under DBS, $m = 4$, $\mu = 1$ and  $\theta = 0.5$.}
\label{fig:ffprime}
\end{figure}

\begin{conjecture}
\label{con:fkci} Under either DBS or DES$,$ 
$\ell(\lambda)$ is increasing and strictly convex in $\lambda,$ having a smooth  derivative satisfying 
  $\ell'(\lambda) \nearrow 1$ as $\lambda \nearrow \infty.$
\end{conjecture}

In the sequel it is assumed that Conjecture \ref{con:fkci} holds, in which case  
  (\ref{eq:nlp}) would be a nonlinear optimization problem 
with strictly convex objective and linear constraints. 
Hence, by standard results it would have a unique global optimum  $\boldsymbol{\lambda}^*$  determined by the above KKT conditions.

We assume henceforth,  by reordering if necessary, that $\alpha_1 \leqslant \cdots \leqslant \alpha_n$. Note that under DES this is equivalent to ordering basic nodes in nonincreasing order of server speed, so $\mu_1 \geqslant \cdots \geqslant \mu_n$. 

Under Conjecture \ref{con:fkci}, it follows that for any basic node $k$ and $\alpha \in (\alpha_k, 1)$ 
the equation $\ell_k'(\lambda_k) = \alpha$ in variable $\lambda_k$ has a unique root $\lambda_k^*(\alpha) > 0$, which is a smooth increasing function of $\alpha$ satisfying $\lambda_k^*(\alpha) \searrow 0$ as $\alpha \searrow \alpha_k$  and $\lambda_k^*(\alpha) \nearrow \infty$ as $\alpha \nearrow 1$.
We shall find it convenient to extend the domain of such a function from $(\alpha_k, 1)$ to 
$\mathbb{R}_+ \triangleq [0, \infty)$ by further defining
\begin{equation}
\label{eq:lambdakstext}
\lambda_k^*(\alpha)
 \triangleq  
\begin{cases}
0 & \quad \textup{for } 0 \leqslant \alpha \leqslant \alpha_k \\
\infty & \quad \textup{for } \alpha \geqslant 1.
\end{cases}
\end{equation}
We shall write, for $\alpha \in \mathbb{R}_+$, 
$\boldsymbol{\lambda}^*(\alpha) \triangleq (\lambda_k^*(\alpha))_{k=0}^n$, with 
$\lambda_0^*(\alpha) \triangleq \lambda - \sum_{k=1}^n \lambda_k^*(\alpha)$.

Let $\Lambda^*(\alpha) \triangleq
\sum_{k=1}^n \lambda_k^*(\alpha)$, and 
note that
\begin{equation}
\label{eq:Lambdast}
\Lambda^*(\alpha) = 
\begin{cases}
0 & \textup{if } 0 \leqslant \alpha \leqslant \alpha_1 \\
\displaystyle \sum_{k=1}^l \lambda_k^*(\alpha) & \textup{if } \alpha_l \leqslant \alpha \leqslant \alpha_{l+1}, \quad l = 1, \ldots, n-1 \\
\displaystyle \sum_{k=1}^n \lambda_k^*(\alpha) & \textup{if } \alpha \geqslant \alpha_{n}.
\end{cases}
\end{equation}

Fig.\ \ref{fig:bigLambdastar} displays the function $\Lambda^*(\alpha)$ for a given instance in the DBS and  DES cases.
The equations $\ell_k'(\lambda_k) = \alpha$ and (\ref{eq:Rstar}) have been solved using the MATLAB \texttt{fzero} function.
Note that $\Lambda^*(\alpha)$ is an increasing function of $\alpha$ that is smooth in the DBS case and piecewise smooth with $n+1$ pieces in the DES case, consistently with (\ref{eq:Lambdast}), with $\Lambda^*(\alpha) \nearrow \infty$ as $\alpha \nearrow 1$.
 
\begin{figure}[!htb]
\centering
\includegraphics[height=1.7in]{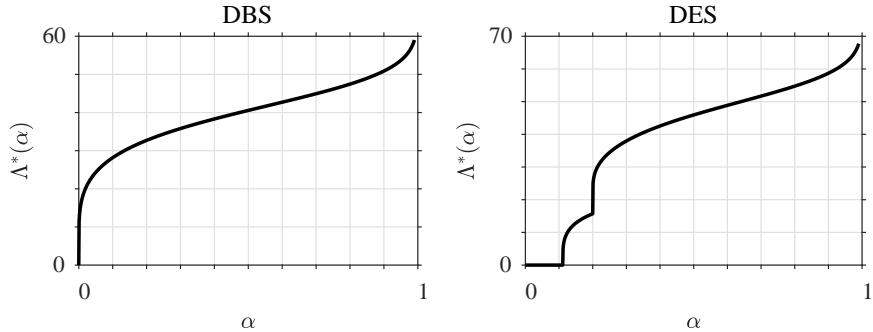}
\caption{$\Lambda^*(\alpha)$ for the instance with 
$n = 2$, $\theta = 0.5$, $(m_1, m_2) = (5, 10)$, $(\mu_1, \mu_2) = (4, 2)$.}
\label{fig:bigLambdastar}
\end{figure}

It follows from the above that, for any $\lambda > 0$,  the equation in $\alpha$
\begin{equation}
\label{eq:Rstar}
\Lambda^*(\alpha) 
= \lambda
\end{equation}
 has a unique root $C^*(\lambda) \in (\alpha_1, 1)$, which can be numerically computed, e.g., 
 by the \emph{bisection method}.
 See \citet[Ch.\ 2.1]{burdenFaires11}.
Note that $C^*(\lambda)$ is an increasing function of $\lambda$ that  is smooth under DBS and piecewise smooth with $n$ pieces under DES, with $C^*(0^+) = \alpha_1$ and $C^*(\lambda) \nearrow 1$ as $\lambda \nearrow \infty$, which follows from corresponding properties of $\Lambda^*(\alpha)$. Let also $C^*(0) \triangleq \alpha_1$.
 
The following result characterizes the optimal BS and identifies relevant properties of it.
Let $\alpha^* \triangleq \min \{C, C^*(\lambda)\}$. 

\begin{proposition}
\label{pro:obsastar}
Assuming Conjecture $\ref{con:fkci},$ the following holds$:$
\begin{itemize}
\item[\textup{(a)}]
$\boldsymbol{\lambda}^*(\alpha^*)$ is the optimal BS$;$
\item[\textup{(b)}] $\lambda_0^*(\alpha^*) > 0$ if and only if $C < C^*(\lambda);$
\item[\textup{(c)}] $\lambda_k^*(\alpha^*) > 0$ if and only if $\alpha^* > \alpha_k,$ for any basic node $k = 1, \ldots, n.$
\end{itemize}
\end{proposition}
\begin{proof}
(a, b) The proofs of the two parts are intertwined.
In the case $C < C^*(\lambda)$ where $\alpha^* = C$, 
 it follows from the above that
\[
\sum_{k=1}^n \lambda_k^*(\alpha^*) < \lambda,
\]
whence $\lambda_0^*(\alpha^*) > 0$.

As for the case $C \geqslant C^*(\lambda)$ where $\alpha^* = C^*(\lambda)$, 
we have
\[
\sum_{k=1}^n \lambda_k^*(\alpha^*) = \lambda,
\]
whence $\lambda_0^*(\alpha^*) = 0$.

Furthermore, it is readily verified that $\boldsymbol{\lambda}^*(\alpha^*)$ and $\alpha^*$ satisfy  conditions (\ref{eq:kktc}), which, assuming Conjecture \ref{con:fkci}, are sufficient for global optimality.

Part (c) follows from the definition of $\lambda_k^*(\alpha)$ and (\ref{eq:lambdakstext}). 
\end{proof}

\begin{remark}
\label{re:obsastar}
The following qualitative properties of the optimal BS are readily derived from Proposition \ref{pro:obsastar}:
\begin{itemize}
\item[\textup{(i)}] Proposition \ref{pro:obsastar}(b) characterizes when the optimal BS resorts to external short-term leased servers  in terms of the threshold cost $C^*(\lambda) \in (\alpha_1, 1)$, as it is optimal to use the external node ($\lambda_0^* > 0$) if and only if the cost of doing so is low enough ($C < C^*(\lambda)$, so $\alpha^* = C$). Hence,  it is not optimal to use the external node if its cost is not lower than the service fee ($C \geqslant 1$). 
Furthermore, for large enough $C$ ($C \geqslant C^*(\lambda)$, so $\alpha^* = C^*(\lambda)$), the optimal BS remains constant, being $\boldsymbol{\lambda}^*(C^*(\lambda))$.
\item[\textup{(ii)}] Since $C^*(0^+) = \alpha_1$ and 
$C^*(\lambda) \nearrow 1$ as $\lambda \nearrow \infty$, for any $C < 1$ the optimal BS uses the external node ($\lambda_0^* > 0$) if and only if the arrival rate $\lambda$ is large enough, viz.\ $\lambda > \lambda^* \triangleq \Lambda^*(C)$, in which case $\alpha^* = C$.
Hence, for $\lambda > \lambda^*$ the $\lambda_k^*$ for basic nodes $k$ remain constant, being equal to $\lambda_k^*(C)$, whereas $\lambda_0^*$ grows linearly in $\lambda$, being $\lambda_0^* = \lambda - \sum_{k=1}^n \lambda_k^*(C)$.
\item[\textup{(iii)}] In light of Proposition \ref{pro:obsastar}(c), we see that in the DBS case where $\alpha_k \equiv 0$ the optimal BS uses all basic nodes ($\lambda_1^*, \ldots, \lambda_n^* > 0$). In contrast, in the DES case where the $\alpha_k$ are positive, the optimal BS may not use some basic nodes: it uses no basic nodes  if and only if  $\alpha^* \leqslant \alpha_1$ (in which case $\alpha^* = C$); 
for $k = 1, \ldots, n-1$, it uses only the $k$ basic nodes with faster servers 
   ($\lambda_l^* > 0$ for $l = 1, \ldots, k$ and 
   $\lambda_l^* = 0$ for $l > k$) if and only if  $\alpha_k < \alpha^* \leqslant \alpha_{k+1}$; 
   and it uses all basic nodes if and only if  $\alpha^* > \alpha_n$.
\item[\textup{(iv)}] The behavior of the optimal BS in the DES case stated in (iii) is more intuitively reformulated in terms of the arrival rate $\lambda$ and the external node usage cost $C$ as follows. 
If $C \leqslant \alpha_{1}$, the optimal BS sends all jobs to the external node.
If $\alpha_{l} < C \leqslant \alpha_{l+1}$ for some $l = 1, \ldots, n-1$, then  the optimal BS does not use the $n-l$ slower basic nodes ($\lambda_{l+1}^* = \cdots = \lambda_n^* = 0$); it uses at least the $k$ faster basic nodes ($\lambda_1^*, \ldots, \lambda_k^* > 0$) for any given $k \leqslant l$ if and only if 
$\lambda > \Lambda^*(\alpha_k)$;  and it uses the external node if and only if $\lambda > \Lambda^*(C)$, in which case it uses the $l$ faster basic nodes.
If $C > \alpha_n$, the optimal BS uses at least the $k$ faster basic nodes  for any given $k$ if and only if 
$\lambda > \Lambda^*(\alpha_k)$;  and it uses the external node if and only if $\lambda > \Lambda^*(C)$, in which case it uses all basic nodes.
\end{itemize}
\end{remark}

\subsection{Computing the optimal BS in large-scale models}
\label{s:cobslsm}
This section discusses computation of the optimal BS in case of a large number of basic nodes, or with some basic nodes having large server pools.

We first consider how computation of the optimal BS scales with $n$.
First, equation (\ref{eq:Rstar}) needs to be solved to obtain (approximately) 
$C^*(\lambda)$ and hence $\alpha^*$.
In light of the definition of $\Lambda^*(\alpha)$ in (\ref{eq:Lambdast}) in terms of vector $\boldsymbol{\lambda}^*(\alpha)$, its evaluation can involve up to 
$n$ evaluations of functions $\lambda_k^*(\alpha)$ corresponding to all basic nodes $k$, so the complexity of this step scales linearly with $n$.
Second, once $C^*(\lambda)$ and $\alpha^*$ are available, the optimal BS is immediately obtained by 
Proposition \ref{pro:obsastar}(a) as $\boldsymbol{\lambda}^*(\alpha^*)$.
Thus, the complexity of approximately computing the optimal BS scales linearly with $n$.

We next consider how having a basic node $k$ with a large number $m_k$ of servers affects computation of the optimal BS.
Note that
the equation $\ell_k'(\lambda_k) = \alpha$ in variable $\lambda_k$ needs to be solved for different values of $\alpha$ to obtain $\lambda_k^*(\alpha)$.
As shown in \ref{s:evalell}, computation of $\ell_k(\lambda_k)$ and $\ell_k'(\lambda_k)$ scales linearly with $m_k$, as it involves the recursive calculation of Erlang-B blocking probabilities for the queue 
M/M/$m_k$/$m_k$ queue.
Hence, other things being equal, computing the optimal BS scales linearly with $m_k$ for a given node $k$.

\subsection{Dependence of optimal BS on model parameters: examples and insights}
\label{s:obsei}
This section explores through examples with $n = 3$ basic nodes
how the optimal BS depends on various model parameters.
Note that the base instances referred to below are those specified in Table \ref{t:bi}.

\begin{figure}[!htb]
\centering
\includegraphics[height=1.7in]{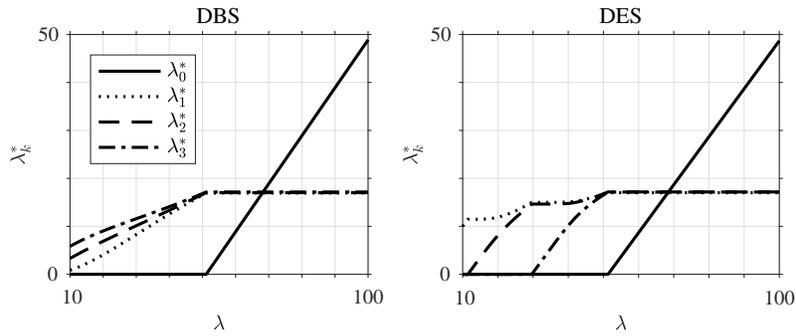}
\caption{Optimal BS vs.\ $\lambda$ for  base instance $1$.}
\label{fig:obsdbse1}
\end{figure}

\subsubsection{Dependence on the arrival rate $\lambda$}
\label{s:obslambda}
We start by investigating the effect of varying the arrival rate $\lambda$.
\paragraph{Base instance $1$}
Fig.\ \ref{fig:obsdbse1} displays the optimal BS vs.\ $\lambda$   for  base instance $1$.
The rate $\lambda_0^*$ of jobs sent to the external node is positive only for large enough $\lambda$ ($\lambda > \lambda^*$),  increasing linearly for $\lambda > \lambda^*$, consistently with Remark \ref{re:obsastar}(ii).

As for the rates $\lambda_k^*$ for basic nodes, their behavior differs in the DBS and DES cases.
In the former, displayed in the left pane of Fig.\ \ref{fig:obsdbse1}, we have
$\lambda_1^* < \lambda_2^* < \lambda_3^*$ for $\lambda < \lambda^*$, while $\lambda_1^* \approx \lambda_2^* \approx \lambda_3^*$ for $\lambda > \lambda^*$.
Thus, when the load is light a larger share of jobs is routed to basic nodes with more servers.

In the DES case, displayed in the right pane of Fig.\ \ref{fig:obsdbse1}, for light loads a larger share of jobs is routed to basic nodes with faster servers. 
Further, nodes $2$ and $3$ are only used for arrival rates $\lambda$ larger than certain critical levels, consistently with Remark \ref{re:obsastar}(iv) as $C > \alpha_3$ in this instance. The behavior for $\lambda > \lambda^*$ is similar to that in the DBS case.

\begin{figure}[!htb]
\centering
\includegraphics[height=1.7in]{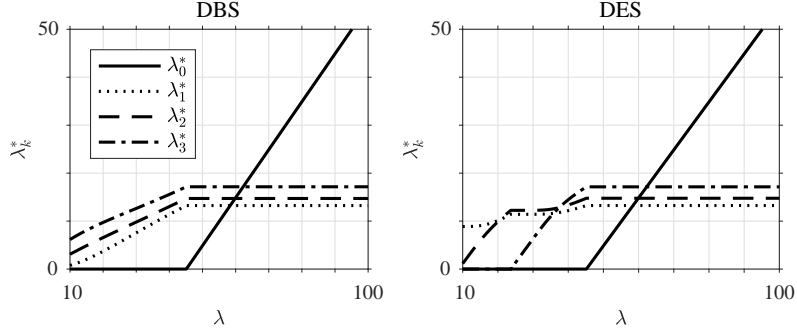}
\caption{Optimal BS vs.\ $\lambda$ for  base instance $2$.}
\label{fig:obsdbse2}
\end{figure}

\paragraph{Base instance $2$}
Fig.\ \ref{fig:obsdbse2} displays the optimal BS vs.\ $\lambda$ for base instance $2$.
Note that 
$\lambda_1^* < \lambda_2^* < \lambda_3^*$ for $\lambda > \lambda^*$, both in the DBS and DES cases, an ordering consistent with total processing capacities, as $m_1 \mu_1 < m_2 \mu_2 < m_3 \mu_3$.
The same holds in the DBS case when the load is light ($\lambda \leqslant \lambda^*$). 
In the DES case we see for $\lambda \leqslant \lambda^*$ a behavior consistent with Remark  \ref{re:obsastar}(iv).

\begin{figure}[!htb]
\centering
\includegraphics[height=1.7in]{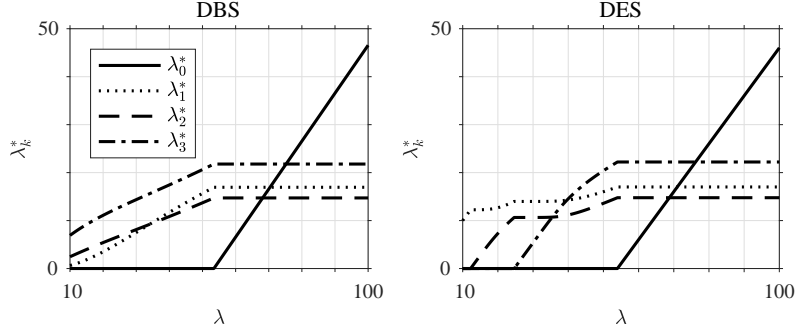}
\caption{Optimal BS vs.\ $\lambda$ for base instance $3$.}
\label{fig:obsdbse3}
\end{figure}

\paragraph{Base instance $3$}
Fig.\ \ref{fig:obsdbse3} displays the optimal BS vs.\ $\lambda$ for base instance $3$.
We see that 
$\lambda_2^* < \lambda_1^* < \lambda_3^*$ for $\lambda > \lambda^*$, both in the DBS and DES cases, again  consistently with total processing capacities, as $m_2 \mu_2 <  m_1 \mu_1 <  m_3 \mu_3$.

When the load is low, we observe in the DBS case the ordering $\lambda_1^* < \lambda_2^* < \lambda_3^*$, consistent with node server pool sizes, as $m_1 < m_2 < m_3$.  
In the DES case we see for low enough $\lambda$ a behavior consistent with Remark  \ref{re:obsastar}(iv).

\paragraph{Insights}
The above examples  suggest the following insights: (1) for high enough load, 
basic nodes with larger total processing capacities $m_k \mu_k$ receive a larger share of traffic, both in the DBS and DES cases;
(2) in the DBS case, for low enough arrival rates basic nodes with larger server pools receive a larger share of traffic;
(3) in the DES case, for low enough arrival rates basic nodes with faster servers  receive a larger share of traffic.

\begin{figure}[!htb]
\centering
\includegraphics[height=1.7in]{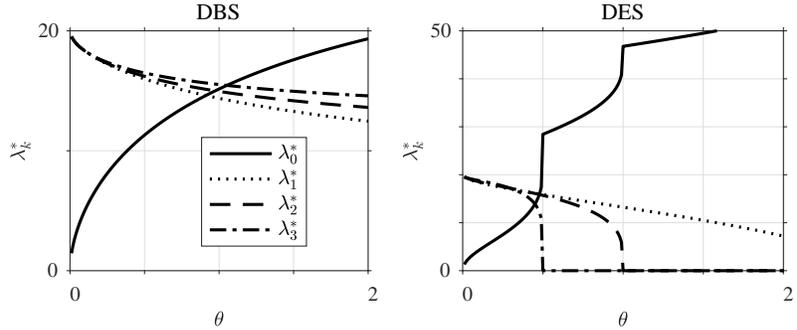}
\caption{Optimal BS vs.\ $\theta$ for base instance $1$.}
\label{fig:obsdbse4}
\end{figure}

\subsubsection{Dependence on the abandonment rate $\theta$}
\label{s:obstheta}
We continue by investigating the effect on the optimal BS of modifying the abandonment rate $\theta$.

\paragraph{Base instance $1$}
Fig.\ \ref{fig:obsdbse4} displays the optimal BS vs.\ $\theta$ as the latter ranges over  $(0, 2]$ for base instance $1$.  
The rate $\lambda_0^*$ of jobs sent to the external node grows with $\theta$, as a smooth concave function in the DBS case, and as a piecewise smooth function with $3$ pieces in the DES case. 

As for the $\lambda_k^*$ for basic nodes, they consequently decrease as $\theta$ grows. Otherwise, their behavior differs substantially in the DBS and DES cases.
In the former, such  $\lambda_k^*$ are very close for very small $\theta$, and get further apart as $\theta$ grows, while maintaining the ordering 
$\lambda_1^* < \lambda_2^* < \lambda_3^*$.
In the latter, such  $\lambda_k^*$ are very close only for small $\theta$. 
For larger $\theta$, they are ordered as $\lambda_3^* < \lambda_2^* < \lambda_1^*$, consistently with server speeds $\mu_3 < \mu_2 < \mu_1$, with   
$\lambda_3^*$ and $\lambda_2^*$ dropping to $0$ at about $\theta \approx 0.5$ and $\theta \approx 1$, respectively.

\begin{figure}[!htb]
\centering
\includegraphics[height=1.7in]{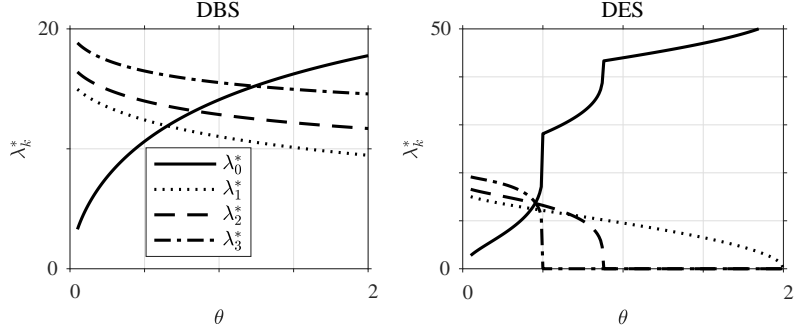}
\caption{Optimal BS vs.\ $\theta$ for base instance $2$.}
\label{fig:obsdbse5}
\end{figure}

\paragraph{Base instance $2$}
Fig.\ \ref{fig:obsdbse5} displays the optimal BS vs.\ $\theta$ for base instance $2$.  
The behavior of $\lambda_0^*$ is similar to that in the previous base instance.
As for the $\lambda_k^*$ for basic nodes, their behavior is also similar, except for the following: (1) in the DBS case, 
such  $\lambda_k^*$ are not very close, while they are also ordered as $\lambda_1^* < \lambda_2^* < \lambda_3^*$, consistently with capacities 
$m_1 \mu_1 < m_2 \mu_2 < m_3 \mu_3$;
and (2) 
in the DES case, such  $\lambda_k^*$ are also ordered as $\lambda_1^* < \lambda_2^* < \lambda_3^*$ for very small $\theta$.

\begin{figure}[!htb]
\centering
\includegraphics[height=1.7in]{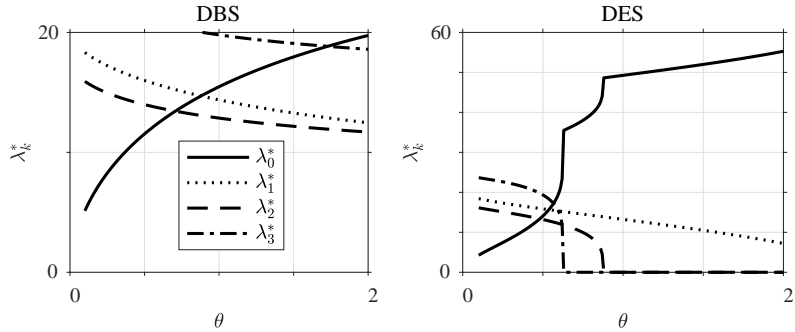}
\caption{Optimal BS vs.\ $\theta$ for base instance $3$.}
\label{fig:obsdbse6}
\end{figure}

\paragraph{Base instance $3$}
Fig.\ \ref{fig:obsdbse6} displays the optimal BS vs.\ $\theta$ for base instance $3$.  
Regarding the $\lambda_k^*$ for basic nodes $k$, we observe that, in the DBS case, 
they are ordered as $\lambda_2^* < \lambda_1^* < \lambda_3^*$, consistently with node capacities 
$m_2 \mu_2 < m_1 \mu_1 < m_3 \mu_3$.
In the DES case, such an ordering also holds for very small $\theta$.

\paragraph{Insights}
The above examples suggest the following insights: (1) in the DBS case, the rate $\lambda_0^*$ of jobs outsourced to the external node increases with $\theta$ as a smooth concave function, while the $\lambda_k^*$ for basic nodes decrease as smooth convex functions, with nodes having a larger total processing capacity receiving a larger share of traffic;
and (2) in the DES case, $\lambda_0^*$ increases steeply with $\theta$, with the $\lambda_k^*$ for basic nodes steeply decreasing; faster nodes receive a larger share of traffic for larger $\theta$, to the point that nodes with slower servers are not used for high enough $\theta$.

Thus, the dependence of the optimal BS on $\theta$ is substantially more pronounced in the DES case than in the DBS case.

\begin{figure}[!htb]
\centering
\includegraphics[height=1.7in]{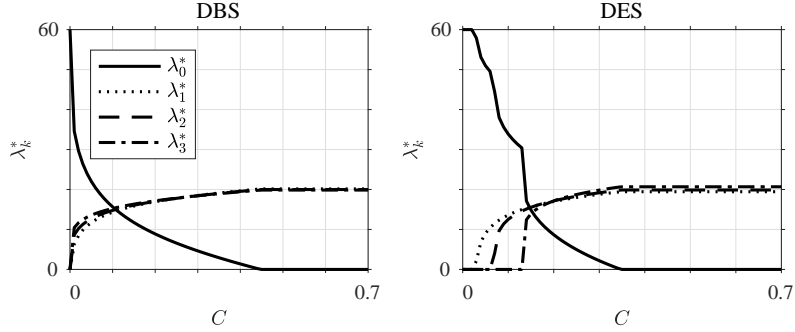}
\caption{Optimal BS vs.\ $C$ for base instance $1$.}
\label{fig:obsdbse7}
\end{figure}

\subsubsection{Dependence on the external node's expected usage cost $C$}
\label{s:obsC}
Finally, we consider the effect on the optimal BS of modifying $C$.

\paragraph{Base instance $1$}
Fig.\ \ref{fig:obsdbse7} displays the optimal BS vs.\ $C$ for base instance $1$.  
The rate $\lambda_0^*$ of jobs outsourced to the external node decreases with $C$, as a convex smooth function in the DBS case for low enough $C$, and as a piecewise smooth function in the DES case, dropping to $0$ for large enough $C$. Such behavior is
consistent with Remark \ref{re:obsastar}(i).

As for the $\lambda_k^*$ for basic nodes, they consequently increase as $C$ grows, remaining constant for large enough $C$.
In the DBS case, such  $\lambda_k^*$ are very close, while appearing 
ordered as $\lambda_1^* < \lambda_2^* < \lambda_3^*$ for very small $C$.
In the DES case, such  $\lambda_k^*$ are very close only for large enough $C$. 
For small $C$, they are ordered as $\lambda_3^* < \lambda_2^* < \lambda_1^*$, consistently with server speeds $\mu_3 < \mu_2 < \mu_1$, and in agreement with Remark \ref{re:obsastar}(iv).

\begin{figure}[!htb]
\centering
\includegraphics[height=1.7in]{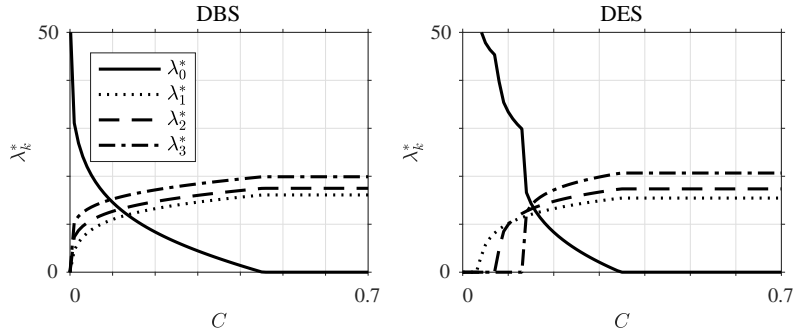}
\caption{Optimal BS vs.\ $C$ for base instance $2$.}
\label{fig:obsdbse8}
\end{figure}

\paragraph{Base instance $2$}
Fig.\ \ref{fig:obsdbse8} displays the optimal BS vs.\ $C$ for base instance $2$.  
The behavior of $\lambda_0^*$ is similar to the previous instance.
As for the $\lambda_k^*$ for basic nodes, 
in the DBS case they are clearly separated and 
ordered as $\lambda_1^* < \lambda_2^* < \lambda_3^*$, as total node capacities.
In the DES case, the same ordering is observed for high enough $C$, while the ordering for small $C$ is consistent with node server speeds, as indicated in Remark \ref{re:obsastar}(iv).

\begin{figure}[!htb]
\centering
\includegraphics[height=1.7in]{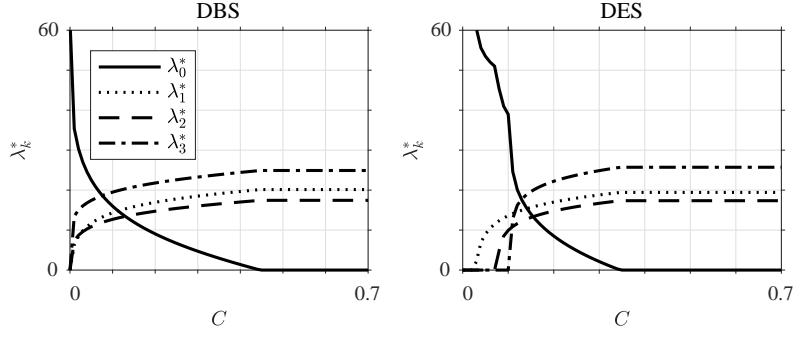}
\caption{Optimal BS vs.\ $C$ for base instance $3$.}
\label{fig:obsdbse9}
\end{figure}

\paragraph{Base instance $3$}
Fig.\ \ref{fig:obsdbse9} displays the optimal BS vs.\ $C$ for base instance $3$.  
The behavior of $\lambda_0^*$ is similar to that in the previous instances.
As for the $\lambda_k^*$ for basic nodes, 
in the DBS case they are clearly separated and 
ordered as $\lambda_2^* < \lambda_1^* < \lambda_3^*$, as total node capacities.
In the DES case, the same ordering is observed for high enough $C$, while the ordering for small $C$ is consistent with nodes server speeds, as in Remark \ref{re:obsastar}(iv).

\paragraph{Insights}
The above examples suggest the following insights: (1) in the DBS case, the rate $\lambda_0^*$ of jobs outsourced to the external node decreases with $C$ as a smooth convex function for $C < C^*(\lambda)$, while the $\lambda_k^*$ for basic nodes increase as smooth concave functions for such $C$, with nodes having a larger total processing capacity receiving a larger share of traffic;
and (2) in the DES case, $\lambda_0^*$ decreases with $C$ as a piecewise-smooth convex function with $n+1$ pieces for $C < C^*(\lambda)$; the $\lambda_k^*$ for basic nodes vanish for small enough $C$, and then increase as smooth concave functions; faster nodes receive a larger share of traffic for smaller $C$, while for larger $C$ the ordering is as under DBS.

\section{Index policies}
\label{s:ip}
\subsection{IO index policy}
\label{s:ioip}
This section  derives an index policy for problem (\ref{eq:panew}) based on 
 indices $\varphi_k^{\scriptscriptstyle \mathrm{IO}}(i_k)$ that measure the probability that an individual job sent to basic node $k$ when this lies in state $i_k$ abandons due to missing its deadline.
The resulting \emph{IO index policy} makes dynamic resource allocation and routing decisions for a job accounting only for their impact on itself. 
The analysis below considers a basic node in isolation, dropping the label $k$ from the notation.
 
Thus, consider an M/M/$m+$M queue with arrival rate $\lambda$, service rate $\mu$ per server, and  abandonment rate $\theta$. 
Denote by $p_{\mathrm{ab}}(i)$ the probability that an arriving job finding $i$ jobs present abandons due to missing its deadline.
The \emph{IO index} is $\varphi^{\scriptscriptstyle  \mathrm{IO}}(i) \triangleq p_{\mathrm{ab}}(i)$, which is evaluated below.
Recall that $L(i)$ is the total abandonment rate in state $i$, and let $D(i)$ be the \emph{total death rate} in state $i$, so $L(i) \triangleq (i-m)^+ \theta$ under DBS and $L(i) \triangleq i \theta$ under DES, and 
\begin{equation}
\label{eq:barmu}
D(i) \triangleq \min(i, m) +L(i).
\end{equation}

\begin{proposition}
\label{pro:ioi} Under either DBS or DES$,$
\begin{equation}
\label{eq:ioindxue}
p_{\mathrm{ab}}(i) = \frac{L(i+1)}{D(i+1)}, \quad i \in \mathbb{Z}_+.
\end{equation}
\end{proposition} 
\begin{proof}
Start with the 
DBS case. 
Consider an arrival that finds $i$ jobs in system. 
To simplify the evaluation of $p_{\mathrm{ab}}(i)$ we assume that there are no further arrivals, which is without loss of generality given the assumption that jobs are scheduled in FCFS order (see \S\ref{s:md}).
If the arrival finds a free server it does not abandon, so
 $p_{\mathrm{ab}}(i) = 0 = L(i+1)/D(i+1)$ for $i < m$.
If $i = m$, so the arrival finds all servers busy and no queue, 
if will abandon if and only if its deadline expires before any service is completed, which happens with probability (w.p.) $\theta / (\theta + m \mu) = L(m+1)/D(m+1)$, so 
 (\ref{eq:ioindxue}) also holds for $i = m$. 
If $i = m+j$, so the arrival finds all servers busy and $j \geqslant 1$ jobs waiting, the next event will be one of the following three: 
the arriving job abandons, w.p.\ $\theta/{D(m+j+1)}$; or
a service is completed, w.p.\ $m \mu/{D(m+j+1)}$; or
one of the $j$ waiting jobs found on arrival abandons, w.p.\ $j \theta/{D(m+j+1)}$.
By conditioning on the next event and using such probabilities one readily obtains the recursion
\[
p_{\mathrm{ab}}(m + j) = 
\frac{\theta}{D(m + j+1)} + \frac{j \theta + m \mu}{D(m + j+1)} p_{\mathrm{ab}}(m + j-1), \quad j = 1, 2, \ldots,
\]
whose solution (taking into account the value obtained for $p_{\mathrm{ab}}(m)$), is 
\[
p_{\mathrm{ab}}(m + j) = \frac{(j+1) \theta}{(j+1) \theta + m \mu} =
\frac{L(m + j+1)}{D(m + j+1)}, \quad j = 0, 1, \ldots.
\]

In the DES case, we have, for $i < m$, that $p_{\mathrm{ab}}(i)$ is the probability that the arriving job, which enters service immediately, completes service before its deadline expires, so 
 $p_{\mathrm{ab}}(i) = \theta/(\theta + \mu) = i \theta/ (i \theta + i \mu) = L(i+1)/D(i+1)$.
Further, arguing along the same lines as in the DBS case yields  the recursion
\[
p_{\mathrm{ab}}(m + j) = 
\frac{\theta}{D(m + j+1)} + \frac{(m+j) \theta + m \mu}{D(m + j+1)} p_{\mathrm{ab}}(m + j-1), \quad j = 0, 1, \ldots,
\]
whose solution (taking into account the value obtained for $p_{\mathrm{ab}}(m-1)$), is
\[
p_{\mathrm{ab}}(m + j) = \frac{(m+j+1) \theta}{(m+j+1) \theta + m \mu} = 
\frac{L(m + j+1)}{D(m + j+1)}, \quad j = 0, 1, \ldots.
\]
\end{proof}

\begin{remark}
\label{re:ioip}
\begin{itemize}
\item[\textup{(i)}] It is immediate that, 
in both the DBS and DES cases$,$ the IO index $\varphi^{\scriptscriptstyle  \mathrm{IO}}(i)$ satisfies target properties $P_1$--$P_6$ (see Table \ref{t:indexp} in \S\ref{s:dpeop}). 
\item[\textup{(i)}]
 $P_2$ holds as $\varphi^{\scriptscriptstyle  \mathrm{IO}}(i)$ does not depend on the arrival rate  $\lambda$. 
\item[\textup{(iii)}] $P_5$ holds with 
$\varphi^{\scriptscriptstyle  \mathrm{IO}}(i) \searrow 0$ in the DBS case and $\varphi^{\scriptscriptstyle  \mathrm{IO}}(i) \searrow \theta/(\theta + \mu)$ in the DES case as $m \nearrow \infty$.
\end{itemize}
\end{remark}
 
\subsection{PI index policy}
\label{s:model2bs}
We next turn to the PI method (see, e.g., \cite{krish90}) which has not been applied before to the present model. 
This method involves two stages: (1) finding the optimal BS of the arrival stream, which is addressed in \S\ref{s:obsm}; and (2) carrying out one step of the PI algorithm for MDPs starting from the optimal BS, which is guaranteed to produce a better policy.

The PI method yields an index policy based on indices $\varphi_k^{\scriptscriptstyle  \mathrm{PI}}(i_k)$ for basic nodes $k = 1, \ldots, n$ that are defined as follows.
Consider the M/M/$m_k$+M queue corresponding to basic node $k$ with arrival rate $\lambda_k^*$, as determined 
by the optimal BS (see \S\ref{s:obsm}). For such a model, let
\[
b_k^*(i) \triangleq \Ex_i\left[\int_0^\infty \{L_k(X_k(t)) - \ell_k(\lambda_k^*)\} \, dt\right], \quad i = 0, 1, 2, \ldots, 
\]
where $\Ex_i[\cdot]$ denotes expectation starting from state $X_k(0) = i$ and $\ell_k(\lambda_k^*)$ is 
the mean abandonment rate as in \S\ref{s:obsm}.
Quantity $b_k^*(i)$, which is known as a \emph{bias} or \emph{relative cost}, measures the expected total incremental number of abandonments when starting from state $i$ relative to those when starting from steady state. Now, the PI index for basic node $k$ is given by 
\begin{equation}
\label{eq:nukpii}
\varphi_k^{\scriptscriptstyle  \mathrm{PI}}(i) \triangleq b_k^*(i+1) - b_k^*(i), \quad i = 0, 1, 2, \ldots, 
\end{equation}

We next address numerical evaluation of the PI index for a basic node, whose label $k$ is dropped from the notation below. 
To evaluate the PI index $\varphi^{\scriptscriptstyle  \mathrm{PI}}(i)$ 
we must solve the \emph{Poisson equations} (see, e.g., \citet{glynnMeyn96}) for the corresponding M/M/$m+$M queue with arrival rate $\lambda^*$ ---as determined by the optimal BS--- and cost rates 
$L(i)$. Such equations are given by 
\begin{equation}
\label{eq:pe1piim}
\begin{split}
\ell(\lambda^*) + \lambda^* b(0)  & = \lambda b(1) \\
\ell(\lambda^*) + \{\lambda^* + D(i)\} b(i)  & = L(i) + \lambda^* b(i+1) + D(i) b(i-1), \quad  i = 1, 2, \ldots,
\end{split}
\end{equation} 
and are immediately reformulated in terms of the PI index $\varphi^{\scriptscriptstyle  \mathrm{PI}}$ in (\ref{eq:nukpii}) as
\begin{equation}
\label{eq:pepiim}
\begin{split}
\ell(\lambda^*) - \lambda^* \varphi^{\scriptscriptstyle  \mathrm{PI}}(0)  & = 0 \\
\ell(\lambda^*) - \lambda^* \varphi^{\scriptscriptstyle  \mathrm{PI}}(i) + D(i)
\varphi^{\scriptscriptstyle  \mathrm{PI}}(i-1) & = L(i), \quad  i = 1, 2, \ldots,
\end{split}
\end{equation}
where $D(i)$ is as in (\ref{eq:barmu}) and $\ell(\lambda^*)$ is the node's mean abandonment rate.

\begin{remark}
\label{re:piip}
\begin{itemize}
\item[\textup{(i)}] Although the first-order linear recursion
 (\ref{eq:pepiim}) cannot be solved in closed form in terms of elementary functions, 
 it gives an efficient means of computing index values
$\varphi^{\scriptscriptstyle  \mathrm{PI}}(0), \ldots, \varphi^{\scriptscriptstyle  \mathrm{PI}}(i)$ in $O(i)$ time, given $\lambda^*$ and $\ell(\lambda^*)$. 
However, such a recursion suffers from numerical instability for large states $i$ (see \citet{nmnetgcoop14}). 
\item[\textup{(ii)}]
If  the optimal BS does not use the node of concern ($\lambda^* = 0$ and hence $\ell(\lambda^*) = 0$), which can happen in the DES case (see Remark \ref{re:obsastar}(iii, iv)), then (\ref{eq:pepiim}) yields
$\varphi^{\scriptscriptstyle  \mathrm{PI}}(i) = L(i+1)/D(i+1)$, 
whence (see Proposition \ref{pro:ioi}) the node's PI index $\varphi^{\scriptscriptstyle  \mathrm{PI}}(i)$ reduces to its IO index $\varphi^{\scriptscriptstyle  \mathrm{IO}}(i)$.
\item[\textup{(iii)}] Observation of the PI index suggests that 
it satisfies target property $P_1$ in Table \ref{t:indexp}. Yet, establishing such a property is a challenging analytical task, as it involves proving corresponding properties of the mean abandonment rate $\ell(\lambda)$ that the author has not found in the literature.
\item[\textup{(iv)}] Being based on the optimal BS, a node's PI index $\varphi^{\scriptscriptstyle  \mathrm{PI}}(i)$ incorporates all model parameters, including server rates and server pool sizes of other nodes as well as the external node's usage cost.
\item[\textup{(v)}] A node's PI index inherits properties of the optimal BS: other things being equal, it does not change for large enough external node's cost $C$, nor for large enough arrival rate $\lambda$. See Remark \ref{re:obsastar}(i, ii).
\end{itemize}
\end{remark}

\subsection{RB index policy}
\label{s:rbip}
The RB policy results by applying to the present model the  
index policy proposed in \citet{whit88b}  for the general MARBP.
To deploy such an approach, 
problem (\ref{eq:panew}) needs to be cast into the framework of the MARBP, which
 concerns the optimal dynamic activation of a collection of stochastic \emph{projects} modeled as RBs ---binary-action (active or passive) MDPs--- subject to given activation constraints. 
The required reformulation was first given in \citet[\S8.1]{nmmp02} in a broader model for optimal control of admission  and routing to parallel queues, and further developed in \citet{nmnetcoop07}.

To make this paper self-contained, the reformulation of the present model as an MARBP and the RB index policy are outlined in \ref{a:rbip}.
  
\subsection{Dependence of the indices on the state and model parameters}
\label{s:civdpni}
This section explores through examples 
how the routing indices considered depend on the node's state and 
various model parameters, and draws insights from the results obtained.
In particular, satisfaction of target properties $P_1$--$P_6$ (see Table \ref{t:indexp} in \S\ref{s:dpeop}) is investigated.
The  base instance numbers referred to below are as in Table \ref{t:bi} in \S\ref{s:obsei}.

\begin{figure}[!htb]
\centering
\includegraphics[height=1.7in]{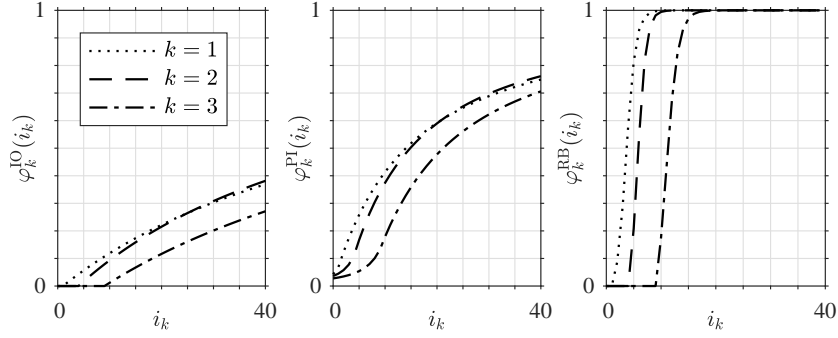}
\caption{Indices vs.\ state for base instance 3 under DBS.}
\label{fig:indbse3}
\end{figure}

\subsubsection{Dependence on the node's state $i_k$}
\label{s:dotsik}
Fig.\ \ref{fig:indbse3} plots the indices for each of the basic nodes $k = 1, 2, 3$ vs.\ the state $i_k$ for base instance 3 under DBS.
The plot is consistent with target property $P_1$ in Table \ref{t:indexp}.
Furthermore, note that the IO indices grow relatively slowly with the state, the PI indices grow faster, and the RB indices grow steeply. In particular,  
the plot shows that such indices satisfy
\begin{equation}
\label{eq:rbsi}
\begin{split}
\varphi_k^{\scriptscriptstyle  \mathrm{IO}}(i) < \varphi_k^{\scriptscriptstyle  \mathrm{PI}}(i) < \varphi_k^{\scriptscriptstyle  \mathrm{RB}}(i), & \quad i \geqslant m_k.
\end{split}
\end{equation}
A similar behavior is observed in Fig.\ \ref{fig:indbse3b} for the same instance under DES.

\begin{figure}[!htb]
\centering
\includegraphics[height=1.7in]{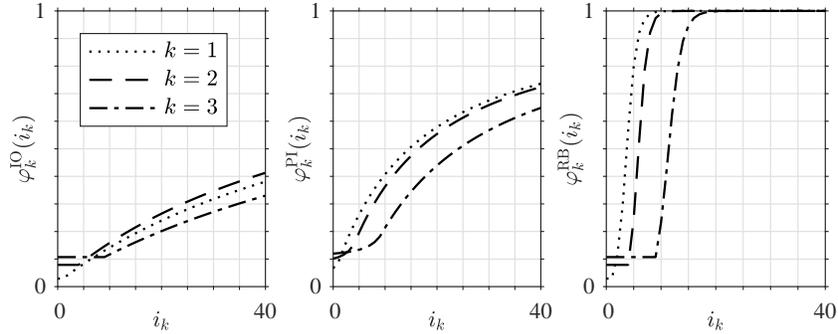}
\caption{Indices vs.\ state for base instance 3 under DES.}
\label{fig:indbse3b}
\end{figure} 

An implication of (\ref{eq:rbsi}), and of the relative index magnitudes as illustrated in the plots, is that the RB policy will tend to outsource more jobs to the external node than  the PI policy, and this more than the IO  policy. 
Note that a key insight from optimal admission control to service  systems (see  \citet{stid85}) is that IO admission policies admit more jobs than is socially optimal. 

\begin{figure}[!htb]
\centering
\includegraphics[height=1.7in]{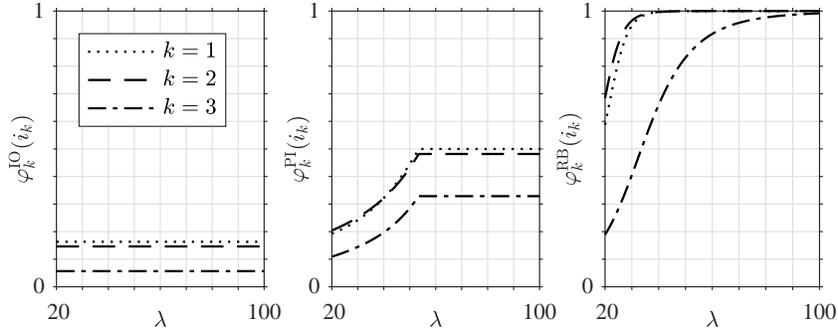}
\caption{Indices for $i_k \equiv 14$ vs.\ $\lambda$ for base instance $3$ under DBS.}
\label{fig:indbsLam3}
\end{figure}

\subsubsection{Dependence on the arrival rate $\lambda$}
\label{s:dotarl}
Fig.\ \ref{fig:indbsLam3} plots the indices for state $i_k = 14$ vs.\ $\lambda$ for base instance $3$ under DBS.
The plot is consistent with property $P_2$ in Table \ref{t:indexp}. 
Note that the IO index remains constant as $\lambda$ varies, as it does not incorporate such a parameter.
As for the PI index, it grows as a convex function up to a critical value, and then remains constant for larger $\lambda$, consistently with
Remark \ref{re:piip}(v).
Regarding the RB index, it grows steeply to one as $\lambda$ grows.
The corresponding DES case is not shown as results are similar.

\begin{figure}[!htb]
\centering
\includegraphics[height=1.7in]{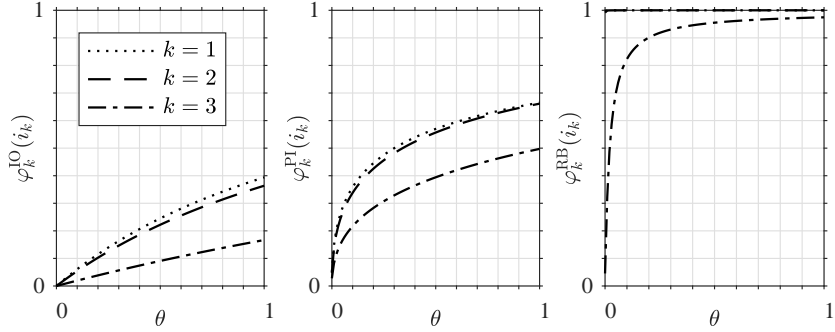}
\caption{Indices for state $\mathbf{i} = (14, 14, 14)$ vs.\ $\theta$ for base instance $3$ under DBS.}
\label{fig:indbsThet3}
\end{figure}

\subsubsection{Dependence on the abandonment rate $\theta$}
\label{s:dotirt}
Fig.\ \ref{fig:indbsThet3} plots the indices for state $i_k = 14$ vs.\ $\theta$ for base instance $3$ under DBS.
The plot is consistent with property $P_3$ in Table \ref{t:indexp}. 
Note that the IO indices grow relatively slowly as the abandonment rate $\theta$ increases, the PI indices grow faster, and the RB indices grow very steeply, to the point that they are very close to $1$ for nodes $1$ and $2$ even for small $\theta$. 

\begin{figure}[!htb]
\centering
\includegraphics[height=1.7in]{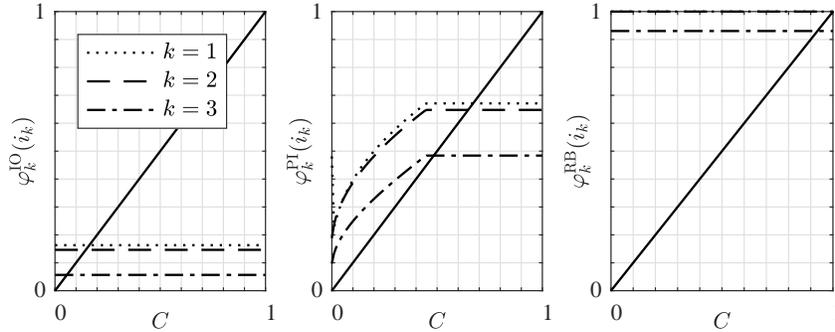}
\caption{Indices for state $\mathbf{i} = (14, 14, 14)$ vs.\ $C$ for base instance $3$ under DBS.}
\label{fig:indbsC3}
\end{figure}

\subsubsection{Dependence on the external node's usage cost $C$}
\label{s:dotesncc}
Fig.\ \ref{fig:indbsC3} plots the indices for state $i_k = 14$ vs.\ $C$ for base instance $3$ under DBS.
The plot is consistent with property $P_4$ in Table \ref{t:indexp}, which can be checked by using the diagonal lines shown. 
Note that both the IO and the RB indices do not vary with $C$, as they do not incorporate such a parameter.
As for the PI indices, they grow as concave functions up to the critical value $C^*(\lambda)$, and then stay constant for larger $C$, consistently with
Remark \ref{re:piip}(v).

\begin{figure}[!htb]
\centering
\includegraphics[height=1.7in]{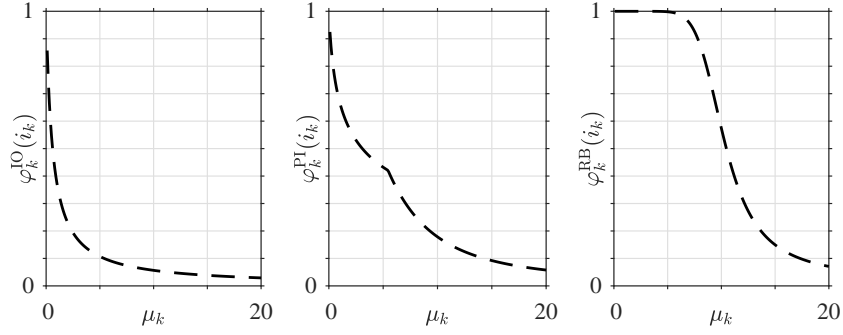}
\caption{Indices for state $i_k =  14$ vs.\ $\mu_k$ for node $k = 2$ in base instance $3$ under DBS.}
\label{fig:indbsmu2}
\end{figure}

\subsubsection{Dependence on the server speed $\mu_k$}
\label{s:dotesr}
Fig.\ \ref{fig:indbsmu2} plots the indices for state $i_k = 14$ vs.\ $\mu_k$ for node $k =2$ of base instance $3$ under DBS.
The plot is consistent with property $P_5$ in Table \ref{t:indexp}. 
Note that the routing priority for node $2$ decreases steeply as its servers get slower (corresponding to  high index values). Such an effect is less pronounced for the PI policy, and still less for the IO policy.

\begin{figure}[!htb]
\centering
\includegraphics[height=1.7in]{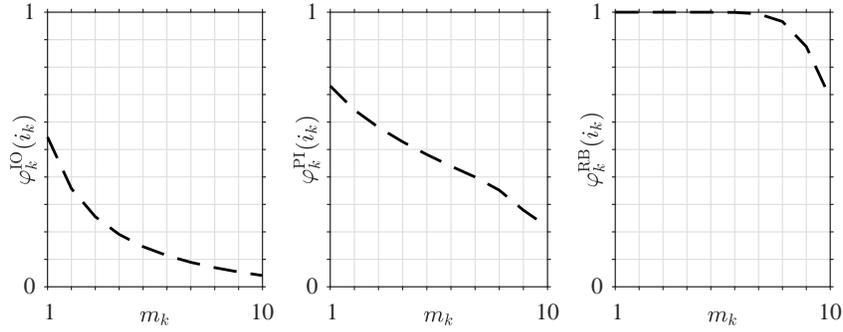}
\caption{Indices for state $i_k =  14$ vs.\ $m_k$ for node $k = 2$ in base instance $3$ under DBS.}
\label{fig:indbsm2}
\end{figure}

\subsubsection{Dependence on the server pool size $m_k$}
\label{s:dotsps}
Fig.\ \ref{fig:indbsm2} plots the indices for state $i_k = 14$ vs.\ $m_k$ for node $k =2$ of base instance $3$ under DBS.
The plot is consistent with property $P_6$ in Table \ref{t:indexp}. 
The routing priority for node $2$ decreases steeply as the number of servers in the node gets smaller (corresponding to high index values). Such an effect is less pronounced for the PI index policy, and even less for the IO policy.

\section{Comparative benchmarking study}
\label{s:cbs}
This section presents a comparative benchmarking study of the four policies considered herein: optimal BS, IO policy, PI policy, and RB policy. 
For such a purpose, a test bed of instances with $n = 2$ basic nodes 
was generated by varying model parameters across a grid of plausible values, while server pool sizes and the rate of node $2$'s servers 
were held constant with $\mathbf{m} = (10, 40)$ and $\mu_2 = 1$.
The other model parameters were varied as follows: 
$\mu_1 \in \{1, 1.5, \ldots, 5\}$, 
$\rho \triangleq \lambda/(m_1 \mu_1 + m_2 \mu_2) \in \{0.9, 1, \ldots, 1.5\}$, 
$\theta \in \{0.2, 0.3, \ldots, 1.1\}$, and
$C \in \{0.1, 0.2, \ldots, 0.8\}$.
Hence, the resulting test bed consists of a total of $5{,}040$ instances, in each the DBS and DES cases.

The performance objective considered was maximization of the
 long-run expected average profit 
per job.
The performance of the optimal policy was used as the benchmark against which the other policies were compared.
For each instance $i$, 
the maximum long-run expected average profit 
per job $z_i^*$ was computed, as well as the value $z_i^{\pi}$ of such an objective under each of the policies $\pi \in \{\mathrm{BS}, \mathrm{IO}, \mathrm{PI}, \mathrm{RB}\}$, truncating the basic nodes' buffer sizes to $80$, as it was checked that increasing them did not change results.
Then, the \emph{relative optimality gap} (percent deviation from optimal performance) 
$100 (z_i^*-z_i^{\pi})/z_i^*$ was evaluated for each policy $\pi$ and instance.

The optimal performance was obtained by solving with CPLEX the linear programming formulation of the DP equations for the truncated model.
The optimal BS's performance was evaluated by solving problem (\ref{eq:nlp}) as discussed in \S\ref{s:obsm}.
As for the performance of the index policies, it was evaluated by solving numerically the corresponding Poisson equations.
The numerical study was based on a MATLAB implementation developed by the author.

\subsection{Overall results}
\label{s:orcomp}
The overall results of the comparative study across all $5{,}040$ instances, both in the DBS and DES cases, are summarized in Tables \ref{t:orog1} and \ref{t:orog2}.

Table \ref{t:orog1} shows the minimum, average and maximum relative optimality gaps for each policy across all instances.
The PI policy tends to outperform the others on such criteria, achieving very low average and maximum relative optimality gaps.
The RB policy closely follows, being also nearly optimal throughout.
The third best policy is the static BS, which is remarkably close to optimal and strongly outperforms the IO policy.

We further observe that the relative optimality gaps of all policies except IO are smaller under DES case than under DBS.
The maximum improvement against the RB policy is small, but it is substantial against the IO policy. 
Yet, both RB and IO can outperform PI, as reflected by the negative minimum improvements. 
As for the improvement of PI over BS, while it is guaranteed to be nonnegative, the
results show that, remarkably, the gains are not major, being at most of $5.33\%$ in the DBS case and $3.94\%$ in the DES case. 

\begin{table}[!htb]
\centering
\begin{tabular}{crrr|}
      & \multicolumn{3}{c}{DBS: optimality gap (\%)} \\ \cline{2-4}
      & \multicolumn{1}{|c}{minimum} & average & maximum \\ \cline{2-4}
BS & \multicolumn{1}{|r}{$0.83$} & $2.99$ &  $5.15$ \\
IO & \multicolumn{1}{|r}{$0.00$} & $6.54$ &  $21.66$  \\
PI & \multicolumn{1}{|r}{$0.00$} & $0.36$ &  $1.40$  \\
RB & \multicolumn{1}{|r}{$0.00$} & $0.44$ &  $2.19$ \\
\cline{2-4}
\end{tabular}
\begin{tabular}{crrr|}
      & \multicolumn{3}{c}{DES: optimality gap (\%)} \\ \cline{2-4}
      & \multicolumn{1}{|c}{minimum} & average & maximum \\ \cline{2-4}
BS & \multicolumn{1}{|r}{$0.00$} & $1.70$ &  $3.79$ \\
IO & \multicolumn{1}{|r}{$0.00$} & $6.57$ &  $22.76$  \\
PI & \multicolumn{1}{|r}{$0.00$} & $0.10$ &  $0.88$  \\
RB & \multicolumn{1}{|r}{$0.00$} & $0.10$ &  $1.36$ \\
\cline{2-4}
\end{tabular}
\caption{Minimum, average and maximum relative optimality gaps of the policies.}
\label{t:orog1}
\end{table} 

Table \ref{t:orog2} shows the minimum, average and maximum percent improvement in the profit per job objective 
of the PI policy vs.\ the other policies. 

\begin{table}[!htb]
\centering
\begin{tabular}{crrr|}
      & \multicolumn{3}{c}{DBS: improvement (\%) of PI} \\ \cline{2-4}
      & \multicolumn{1}{|c}{minimum} & average & maximum \\ \cline{2-4}
BS & \multicolumn{1}{|r}{$0.70$} & $2.72$ &  $5.33$ \\
IO & \multicolumn{1}{|r}{$-0.48$} & $7.02$ &  $27.42$  \\
RB & \multicolumn{1}{|r}{$-1.34$} & $0.09$ &  $2.16$ \\
\cline{2-4}
\end{tabular}
\begin{tabular}{crrr|}
      & \multicolumn{3}{c}{DES: improvement (\%) of PI} \\ \cline{2-4}
      & \multicolumn{1}{|c}{minimum} & average & maximum \\ \cline{2-4}
BS & \multicolumn{1}{|r}{$0.00$} & $1.64$ &  $3.94$ \\
IO & \multicolumn{1}{|r}{$-0.21$} & $7.40$ &  $29.46$  \\
RB & \multicolumn{1}{|r}{$-0.88$} & $0.00$ &  $1.30$ \\
\cline{2-4}
\end{tabular}
\caption{Minimum, average and maximum improvement (\%) of PI vs.\ the other policies.}
\label{t:orog2}
\end{table}

\begin{figure}[!htb]
\centering
\includegraphics[height=2.5in]{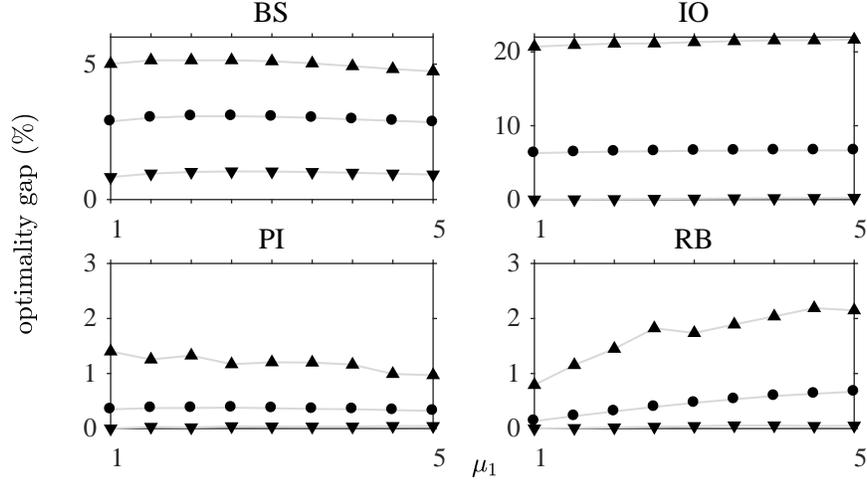}
\caption{Relative optimality gaps vs.\ $\mu_1$ under DBS.}
\label{fig:dbsogapsvmu1}
\end{figure}

\subsection{Effect of changing  model parameters}
\label{s:ecvmp}
To assess the effects of changing model parameters on the policies' relative optimality gaps, the following approach was used: 
for each of the varying parameters and policy, the minimum, average and maximum relative optimality gaps were evaluated across all instances with a given parameter value.

\subsubsection{Effect of changing node $1$'s servers' rate $\mu_1$}
\label{s:en1sr}
Fig.\ \ref{fig:dbsogapsvmu1} displays, for the DBS case,  the smallest, average and largest relative optimality gaps
for the four policies as $\mu_1$ varies, marked as a downward-pointing triangle, a circle, and an upward-pointing triangle, respectively. 

For BS, the gaps remain quite insensitive to changes in $\mu_1$.
The minimum gap ranges between $0.82\%$ and $1.04\%$, 
the average gap between $2.85\%$ and $3.08\%$, and the 
maximum gap between $4.74\%$ and $5.15\%$.

As for IO, the gaps also remain quite stable as $\mu_1$ varies. 
The minimum gap ranges between $0.001\%$ and $0.22\%$, 
the average gap between $6.29\%$ and $6.66\%$, and the 
maximum gap between $20.7\%$ and $21.66\%$. 

Regarding PI, its maximum 
gap slightly improves as $\mu_1$ gets larger.
The minimum gap ranges between $0\%$ and $0.05\%$, 
the average gap between $0.32\%$ and $0.38\%$, and the
maximum gap ranges $0.97\%$ and $1.4\%$. 

The maximum gap of the RB policy worsens as $\mu_1$ grows.
The minimum gap ranges between $0\%$ and $0.06\%$, 
the average gap between $0.14\%$ and $0.67\%$, and the
maximum gap between $0.79\%$ and $2.19\%$. 

Note that for smaller $\mu_1$ RB outperforms PI, whereas for larger $\mu_1$ the opposite holds.
This suggests the appropriateness of the very small routing priority given by the RB index to a node with slow servers, in contrast to the other index policies, as shown in Fig.\ \ref{fig:indbsmu2}. 
It further suggests that the routing priority awarded by the RB policy does not sufficiently increase as the node's servers become faster.

\begin{figure}[!htb]
\centering
\includegraphics[height=2.5in]{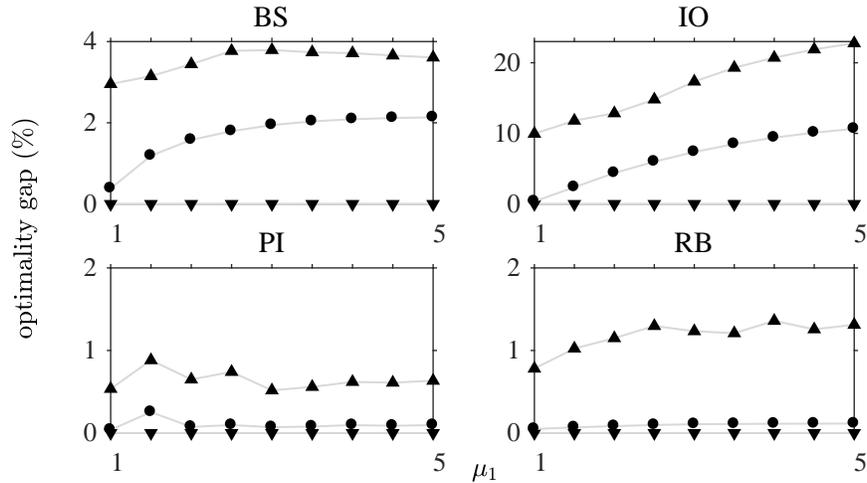}
\caption{Relative optimality gaps vs.\ $\mu_1$ under DES.}
\label{fig:desogapsvmu1}
\end{figure}

Fig.\ \ref{fig:desogapsvmu1} shows corresponding results for the DES case.
The main differences with the results for the DBS case are: 
(1) the gaps are somewhat smaller for all policies; and
(2) the maximum and average gaps for the IO policy increase with $\mu_1$.

\begin{figure}[!htb]
\centering
\includegraphics[height=2.5in]{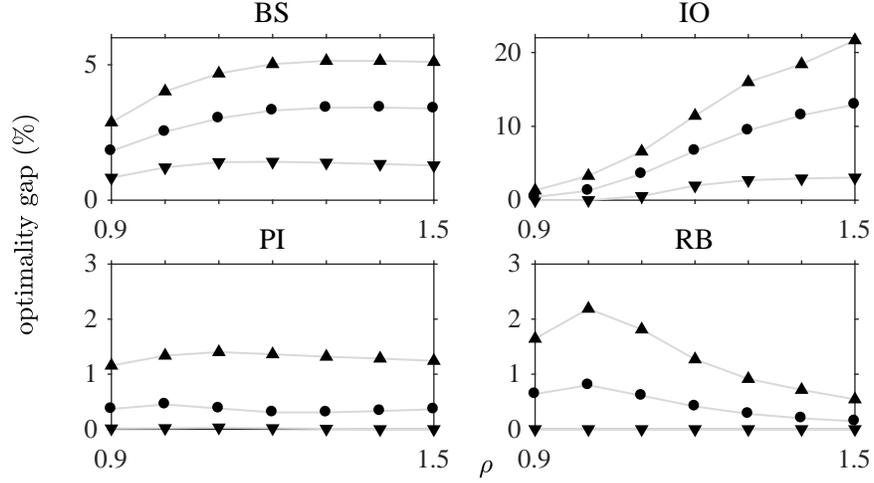}
\caption{Relative optimality gaps vs.\ $\rho$ under DBS.}
\label{fig:dbsogapsvrho}
\end{figure}

\subsubsection{Effect of changing the system's nominal load $\rho$}
\label{s:ecsnlrho}
To assess the effect of changing $\rho$ across the grid of values considered, 
corresponding values of the arrival rate $\lambda$ were generated for each instance 
by taking $\lambda = (m_1 \mu_1 + m_2 \mu_2) \rho $. 
Fig.\ \ref{fig:dbsogapsvrho} displays the results in the DBS case. 

For BS, the three gaps tend to deteriorate as $\rho$ grows.
The minimum gap ranges between $0.83\%$ and $1.42\%$, 
the average gap between $1.82\%$ and $3.42\%$, and the 
maximum gap between $2.87\%$ and $5.15\%$.

As for IO, its gaps deteriorate substantially as $\rho$ grows. 
The minimum gap increases from $0.001\%$ to $3.06\%$, 
the average gap increases from $0.42\%$ to $12.94\%$, and the 
maximum gap increases from $1.32\%$ to $21.66\%$. 
Such behavior is consistent with the fact that such a policy does not incorporate arrival rate information. 

Regarding PI, its  
gaps remain quite stable as $\rho$ varies.
The minimum gap ranges between $0\%$ and $0.03\%$, 
the average gap between $0.3\%$ and $0.45\%$, and the
maximum gap between $1.16\%$ and $1.4\%$. 

As for RB, 
the minimum gap remains at $0\%$, 
the average gap ranges between $0.14\%$ and $0.8\%$, and the
maximum gap between $0.54\%$ and $2.19\%$. 
Note that the maximum and average gaps reach a peak at $\rho = 1$, and then improve as $\rho$ grows.

Note that for smaller $\rho$ PI outperforms RB, whereas for larger $\rho$ RB is better.
This is consistent with the fact that the PI policy is insensitive to the arrival rate as this becomes large, whereas the RB policy steeply decreases the routing priority of basic nodes in heavy traffic, as shown in Fig.\ \ref{fig:indbsLam3}.

\begin{figure}[!htb]
\centering
\includegraphics[height=2.5in]{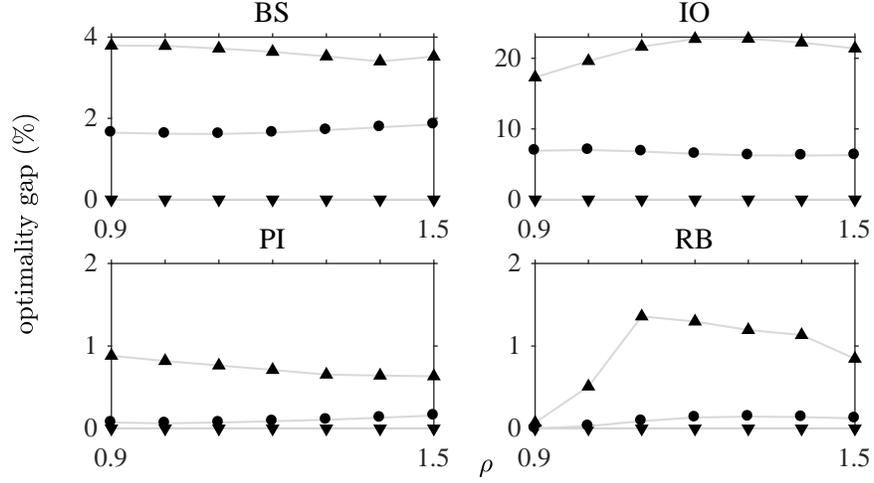}
\caption{Relative optimality gaps vs.\ $\rho$ under DES.}
\label{fig:desogapsvrho}
\end{figure}

Fig.\ \ref{fig:desogapsvrho} shows the results under DES.
The main differences with the DBS case are: 
(1) again, the gaps are somewhat smaller for all policies except for IO; 
(2) the maximum and average gaps for BS and IO do not vary much as $\rho$ changes; and
(3) RB is better than PI for small $\rho$, and PI is better than RB for larger $\rho$ within the range considered.

\begin{figure}[!htb]
\centering
\includegraphics[height=2.5in]{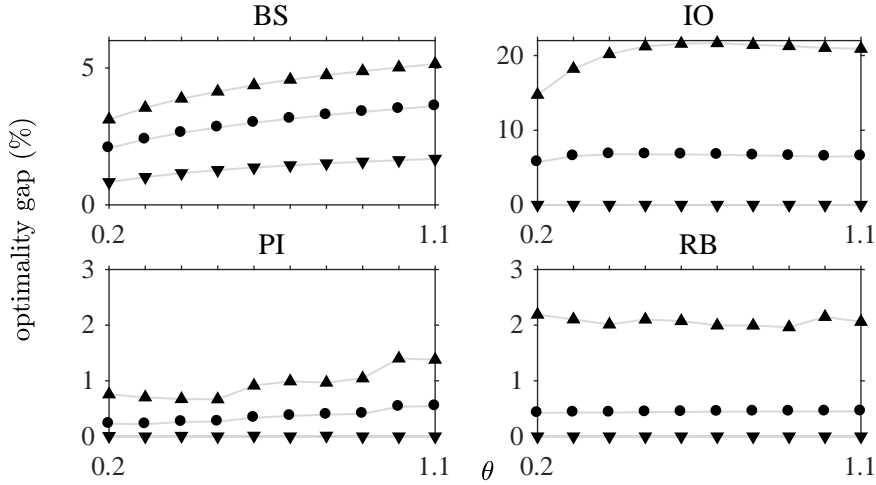}
\caption{Relative optimality gaps vs.\ $\theta$ under DBS.}
\label{fig:dbsogapsvtheta}
\end{figure}

\subsubsection{Effect of changing the abandonment rate $\theta$}
\label{s:ecsnltheta}
Fig.\ \ref{fig:dbsogapsvtheta} displays the smallest, average and largest relative optimality gaps
for each of the policies considered as $\theta$ varies, in the DBS case. 

For BS, the three gaps tend to deteriorate as $\theta$ grows.
The minimum gap increases from $0.83\%$ to $1.68\%$, 
the average gap from $2.07\%$ to $3.6\%$, and the 
maximum gap from $3.12\%$ to $5.15\%$.

As for IO, its worst gap degrades for the smaller values of $\theta$ considered, and then levels off. The average and minimum gaps remain stable. The minimum gap ranges from $0.001\%$ to $0.002\%$, 
the average gap from $5.73\%$ to $6.79\%$, and the 
maximum gap from $14.74\%$ to $21.66\%$. 

Regarding PI, its worst and average gaps degrade slightly as $\theta$ grows.
The minimum gap ranges between $0\%$ and $0.01\%$, 
the average gap between $0.22\%$ and $0.55\%$, and the
maximum gap  between $0.67\%$ and $1.4\%$. 

Concerning RB, its gaps remain stable as $\theta$ varies. 
The minimum gap remains at $0\%$, 
the average gap ranges between $0.14\%$ and $0.8\%$, and the
maximum gap between $1.96\%$ and $2.19\%$. 

Note that PI consistently outperforms RB in the worst case.   

\begin{figure}[!htb]
\centering
\includegraphics[height=2.5in]{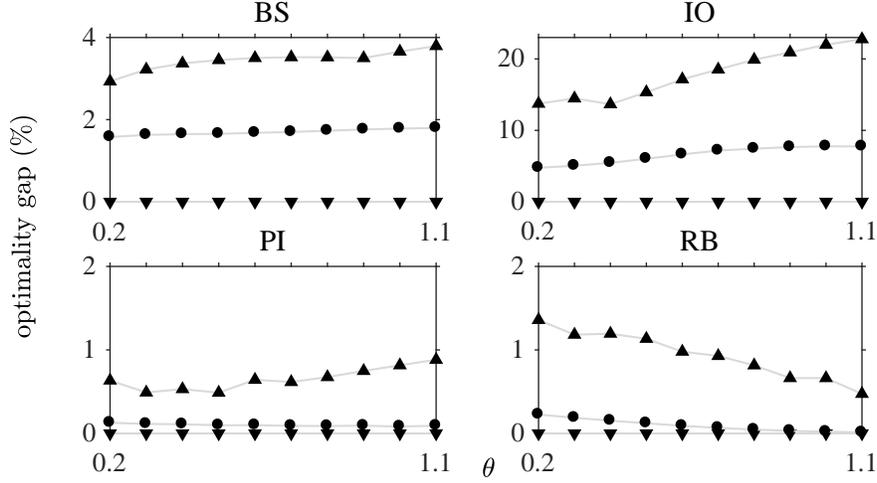}
\caption{Relative optimality gaps vs.\ $\theta$ under DES.}
\label{fig:desogapsvtheta}
\end{figure}

Fig.\ \ref{fig:desogapsvtheta} shows the results under DES.
The main differences under DBS are: 
(1) the gaps are somewhat smaller for all policies except for IO; and
(2) the average gap for BS does not vary much as $\theta$ changes;
(3) the worst gap for IO deteriorates for larger $\theta$; and
(4) the average and worst gaps for RB improves as $\theta$ grows, to the point that RB outperforms PI for larger $\theta$.

\begin{figure}[!htb]
\centering
\includegraphics[height=2.5in]{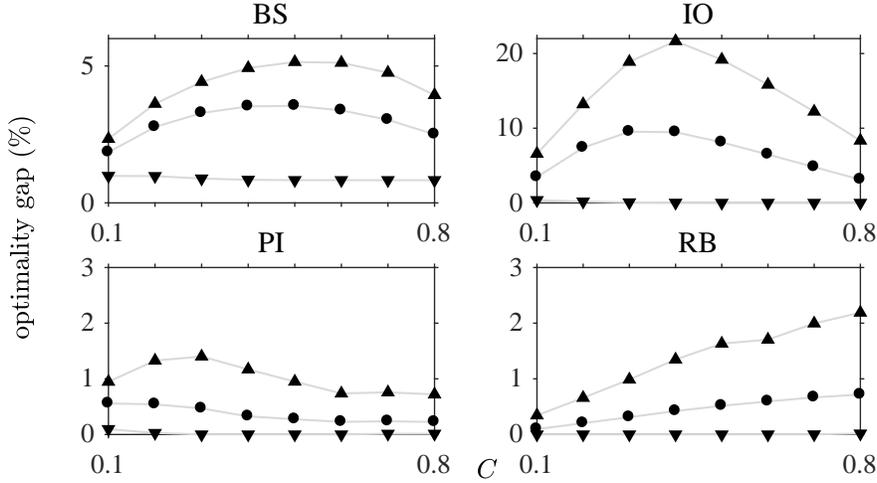}
\caption{Relative optimality gaps vs.\ $C$ under DBS.}
\label{fig:dbsogapsvC}
\end{figure}

\subsubsection{Effect of changing the external node's usage cost $C$}
\label{s:ecesnuc}
Fig.\ \ref{fig:dbsogapsvC} displays the smallest, average and largest relative optimality gaps
for each policy considered as $C$ varies in the DBS case.

For BS, the minimum gap remains stable, ranging from $0.83\%$ to $0.98\%$, 
the average gap increases and then decreases as $C$ grows, 
ranging  from $1.85\%$ to $3.55\%$, and the 
maximum gap has a similar behavior as the average gap, 
ranging from $2.34\%$ to $5.15\%$.

As for IO, its qualitative behavior is similar to that of the BS policy. The minimum gap ranges from $0.001\%$ to $0.37\%$, 
the average gap from $3.12\%$ to $9.53\%$, and the 
maximum gap from $6.58\%$ to $21.66\%$. 

Regarding PI,  its maximum and average gaps improve slightly as $C$ grows.
The minimum gap ranges between $0\%$ and $0.1\%$, 
the average gap decreases from $0.23\%$ to $0.56\%$, and the
maximum gap ranges between $0.72\%$ and $1.4\%$. 

Concerning RB, its maximum and average gaps deteriorate as $C$ grows. 
The minimum gap grows from $0\%$ to $0.02\%$, 
the average gap  from $0.1\%$ to $0.72\%$, and the
maximum gap from $0.34\%$ to $2.19\%$. 

Note that while RB outperforms PI for small $C$, the opposite happens for large $C$.

\begin{figure}[!htb]
\centering
\includegraphics[height=2.5in]{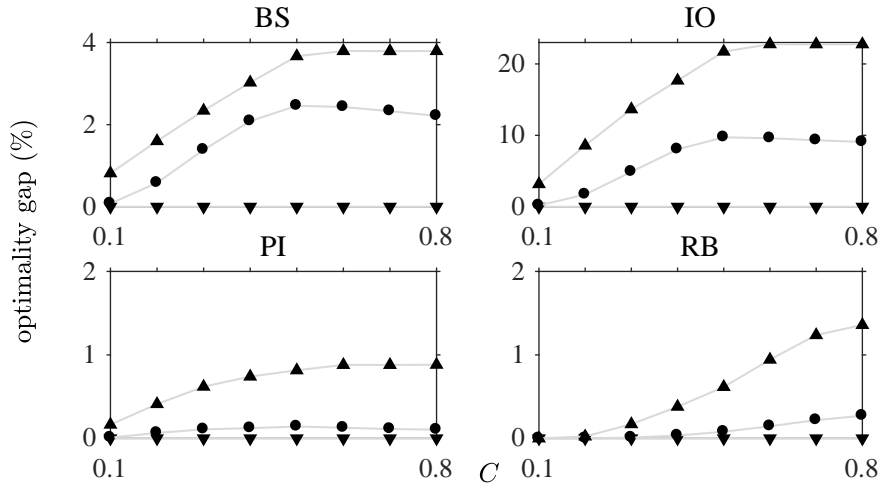}
\caption{Relative optimality gaps vs.\ $C$ under DES.}
\label{fig:desogapsvC}
\end{figure}

Fig.\ \ref{fig:desogapsvC} displays the results in the DES case.
The main differences with the DBS case are that the worst-case gaps deteriorate as $C$ grows for the BS, IO and PI policies, leveling off for large enough $C$. 

\section{Discussion and conclusions}
\label{s:d}
The results presented above allow us to draw some insights about the strengths and weaknesses of the four policies considered.

Concerning the optimal BS policy, it can be efficiently computed, assuming validity of Conjecture \ref{con:fkci}, and its performance is remarkably close to optimal for a static policy, at least for the  range of instances considered.
Its relative optimality gap tends to deteriorate as the abandonment rate grows.

As for the IO policy, it is the easiest to compute, being given in closed form. Yet, despite being a  dynamic policy, its performance is substantially suboptimal, being the worst of the four policies. 
As a result of not incorporating the arrival rate parameter, its relative optimality gap severely degrades as traffic becomes heavier. 
It also fails to appropriately incorporate  other model parameters, as shown in the results in \S\ref{s:ecvmp}. 

Regarding the PI policy, it can be efficiently computed through a simple linear recursion and is overall the best performing of the four, consistently achieving a nearly optimal performance. Yet, it is insensitive to the arrival rate $\lambda$ and to the external node's usage cost $C$ as they become large enough, which causes a slight performance degradation in such cases.

Finally, the RB policy can also be efficiently computed through a linear recursion, and is overall the second-best performing of the four policies, being also consistently near optimal. 
Yet, its index values grow too steeply with the state, which may lead it to outsource too many jobs for large values of $C$, causing a slight deterioration of performance in such cases.

While a theoretical proof that the four routing indices considered satisfy the intuitive target properties listed in Table \ref{t:indexp} remains elusive for PI and RB, 
the above results provide supporting evidence  that such is the case.

This paper raises several challenging issues for future work, such as 
proving (or disproving) Conjecture \ref{con:fkci}, analyzing theoretically whether the 
 policies considered satisfy the target properties in Table \ref{t:indexp}, analyzing their performance, carrying out comparative studies on larger scale instances, and extending the present model to incorporate more complex dynamics.

\section*{Acknowledgements}
This work was supported by the Spanish General Directorate for Scientific and Technical Research (DGICYT)  projects MTM2010-20808 and ECO2015-66593-P. The author has presented
preliminary parts at the 13th IEEE International Conference on High Performance Computing and Communications (HPCC 2011, Banff, Canada), and at the 7th International Conference on Network Games, Control and Optimization (NetGCooP 2014, Trento, Italy), appearing in abridged form in their proceedings (\citet{nmhpcc11,nmnetgcoop14}).

\appendix
\section{Evaluating $\ell(\lambda)$ and $\ell'(\lambda)$}
\label{s:evalell}
\subsection{DBS case}
\label{s:evalelldbs}
This section addresses numerical evaluation of the abandonment rate function $\ell(\lambda) \triangleq \lambda P_{ \mathrm{ab}}(\lambda)$ and its derivative $\ell'(\lambda)$ 
for an M/M/$m+$M queue under DBS, arrival rate $\lambda$, abandonment rate $\theta$, and servers' rate $\mu$, where $P_{ \mathrm{ab}}(\lambda)$ is the abandonment probability.  
This is required for computing the optimal BS by solving problem (\ref{eq:nlp}) in \S\ref{s:obsm}. 
In his seminal analysis of such a model, 
\citet[pp.\ 60--67]{palm57b} gives the formula 
\begin{equation}
\label{eq:apci}
P_{ \mathrm{ab}}(\lambda) = \bigg(\frac{1}{W(\lambda)} + \rho(\lambda) - 1\bigg) \frac{P_{ \mathrm{wait}}(\lambda)}{\rho(\lambda)},
\end{equation}
where 
\begin{equation}
\label{eq:Wldbs}
W(\lambda) \triangleq 1 + \sum_{i=1}^\infty \frac{(\lambda/\theta)^i}{(\beta + 1) \cdots (\beta + i)} = \beta e^{\lambda /\theta} (\lambda /\theta)^{-\beta} \gamma(\beta, \lambda /\theta),
\end{equation}
$\rho(\lambda) \triangleq \lambda / (m \mu)$ is the \emph{nominal load per server}, 
$\beta \triangleq m \mu / \theta$, 
$
\gamma(s, x) \triangleq \int_0^x t^{s-1} e^{-t} \, dt
$
is the \emph{lower incomplete gamma function}, 
and $P_{ \mathrm{wait}}(\lambda)$ is the probability that a random arrival must wait. This is  given by
\begin{equation}
\label{eq:Pwl}
P_{ \mathrm{wait}}(\lambda) = 
\frac{W(\lambda) B_m(r(\lambda))}{1 + \big(W(\lambda)-1\big) B_m(r(\lambda))},
\end{equation}
where $r(\lambda) \triangleq \lambda/\mu$ is the \emph{offered load} and 
$B_{ m}(r)$ is the Erlang-B blocking probability for an M/M/$m$/$m$ queue with offered load $r$. 

For numerical evaluation of $B_{ m}(r)$, one may use the well-known recursion
\begin{equation}
\label{eq:ebf}
B_m(r) = 
\begin{cases}
1 & \textup{ if } m = 0 \\
\displaystyle \frac{r B_{m-1}(r)}{m + r B_{m-1}(r)} & \textup{ if } m \geqslant 1,
\end{cases} 
\end{equation}
or the more numerically stable version
\begin{equation}
\label{eq:ebf2}
B_m^{-1}(r) = 
\begin{cases}
1 & \textup{ if } m = 0 \\
\displaystyle 1 + \frac{m}{r} B_{m-1}^{-1}(r) & \textup{ if } m \geqslant 1.
\end{cases} 
\end{equation}

We next turn to evaluation of the derivative $P_{ \mathrm{ab}}'(\lambda)$, for which we need to evaluate the derivatives  
$B_m'(r)$ and $W'(\lambda)$.
These can be obtained from $B_m(r)$ and $W(\lambda)$ through  the differential relations
\begin{equation}
\label{eq:Bmpr}
B'_m(r) = 
\left(\frac{m}{r} - 1 + B_m(r)\right) B_m(r) 
\end{equation}
and
\begin{equation}
\label{eq:Wpr}
W'(\lambda) = \frac{1 - (1-\rho(\lambda)) W(\lambda)}{\theta \rho(\lambda)}.
\end{equation}

The next result gives a relation useful for evaluating $\ell'(\lambda)$ given $\ell(\lambda)$.
Let 
\begin{equation}
\label{eq:eCf}
C_m(r) \triangleq \frac{B_m(r)}{1-(r/m)[1 - B_m(r)]}.
\end{equation}
Note that, 
if $r < m$,  $C_m(r)$ is the wait probability in an M/M/$m$ queue with offered load $r$, given by Erlang-C formula. Further, (\ref{eq:eCf}) also defines  
$C_m(r)$ for $r \geqslant m$, which satisfies $C_m(r) \to \infty$ as $r \to \infty$.

\begin{proposition}
\label{pro:ellplevdbs}
\[
\ell'(\lambda) = C_m(r(\lambda)) - \Big(\frac{C_m(r(\lambda))}{B_m(r(\lambda))} -1 + (\beta-m) [1-\rho (\lambda )]\Big) \frac{\ell(\lambda)}{\lambda}
 -  \frac{\mu - \theta}{\mu \theta} \frac{\ell^2(\lambda)}{\lambda}.
\]
\end{proposition}
\begin{proof}
The result follows from straightforward algebraic manipulations from 
$\ell'(\lambda) = P_{ \mathrm{ab}}(\lambda) + \lambda  P_{ \mathrm{ab}}'(\lambda)$
using (\ref{eq:apci}), (\ref{eq:Wldbs}), (\ref{eq:Pwl}), (\ref{eq:Bmpr}), 
(\ref{eq:Wpr}) and (\ref{eq:eCf}).
\end{proof}

\subsection{DES case: reduction to the DBS case}
\label{s:evalelldes}
We next turn to the M/M/$m+$M queue under DES. 
\cite{anckGaf62} analyzed the case $m = 1$.
The multi-server case is studied in \citet[\S1.8.4]{gnedKov89}, where
   formulae are given for the steady-state probabilities $p_i$ that there are $i$ jobs in the system. 
Rather than drawing on such results, this section develops an alternative approach, which the author has not found in the literature, exploiting an equivalence between the DBS and DES cases to reduce the calculation of metrics 
$\ell(\lambda)$ and $\ell'(\lambda)$ in the latter case to the former case, addressed  in \ref{s:evalelldbs}.

Denote by $p_i^{{\scriptscriptstyle \mathrm{DBS}}, \mu}(\lambda)$ and $p_i^{{\scriptscriptstyle \mathrm{DES}}, \mu}(\lambda)$
the steady-state probabilities that there are $i$ jobs in system in an M/M/$m+$M queue with parameters $\lambda$,  $\mu$ and $\theta$ under DBS and DES, respectively, where the notation is meant to emphasize the dependence on arrival and server rates. 

\begin{proposition}
\label{pro:pidbspideseq} For $i = 0, 1, \ldots,$
\[
p_i^{{\scriptscriptstyle \mathrm{DES}}, \mu}(\lambda) = p_i^{{\scriptscriptstyle \mathrm{DBS}}, \mu+\theta}(\lambda).
\]
\end{proposition}
\begin{proof}
The result follows immediately from the identity between the 
flow balance equations for the M/M/$m+$M queue under DES and servers' rate $\mu$ and the  M/M/$m+$M queue under DBS and servers' rate $\mu+\theta$.
\end{proof}

The following result, which the author has not found in the literature either, 
 exploits the above equivalence to obtain relations between performance metrics of corresponding DBS and DES systems.
For the M/M/$m+$M queue with parameters $\lambda$,  $\mu$ and $\theta$ under DBS, let 
$M^{{\scriptscriptstyle \mathrm{DBS}}, \mu}(\lambda)$, 
$P_{\mathrm{ab}}^{{\scriptscriptstyle \mathrm{DBS}}, \mu}(\lambda)$ and 
$\ell^{{\scriptscriptstyle \mathrm{DBS}}, \mu}(\lambda)$
 denote, respectively, the mean number of busy servers, the abandonment probability, and the abandonment rate, and let us write as 
 $M^{{\scriptscriptstyle \mathrm{DES}}, \mu}(\lambda)$, 
$P_{\mathrm{ab}}^{{\scriptscriptstyle \mathrm{DES}}, \mu}(\lambda)$ and 
$\ell^{{\scriptscriptstyle \mathrm{DES}}, \mu}(\lambda)$ such metrics in the DES case.

\begin{proposition}
\label{pro:metricsequiv}
\textup{ }
\begin{itemize}
\item[\textup{(a)}] $M^{{\scriptscriptstyle \mathrm{DES}}, \mu}(\lambda) = 
M^{{\scriptscriptstyle \mathrm{DBS}}, \mu+\theta}(\lambda);$
\item[\textup{(b)}] $\displaystyle (\mu+\theta) \big(1-P_{\mathrm{ab}}^{{\scriptscriptstyle \mathrm{DES}}, \mu}(\lambda)\big) = \mu \big(1-P_{\mathrm{ab}}^{{\scriptscriptstyle \mathrm{DBS}}, \mu+\theta}(\lambda)\big);$
\item[\textup{(c)}]
$\displaystyle (\mu + \theta) \ell^{{\scriptscriptstyle \mathrm{DES}}, \mu}(\lambda) = 
\lambda \theta + \mu \ell^{{\scriptscriptstyle \mathrm{DBS}}, \mu+\theta}(\lambda).$
\end{itemize}
\end{proposition}  
\begin{proof}
(a) We have, using Proposition \ref{pro:pidbspideseq}, 
\[
M^{{\scriptscriptstyle \mathrm{DES}}, \mu}(\lambda) = 
\sum_{i=0}^\infty \min(i, m) p_i^{{\scriptscriptstyle \mathrm{DES}}, \mu}(\lambda) =
\sum_{i=0}^\infty \min(i, m) p_i^{{\scriptscriptstyle \mathrm{DBS}}, \mu+\theta}(\lambda) = M^{{\scriptscriptstyle \mathrm{DBS}}, \mu+\theta}(\lambda).
\]

(b) An elementary flow balance argument gives the identities
\[
\lambda \big(1-P_{\mathrm{ab}}^{{\scriptscriptstyle \mathrm{DES}}, \mu}(\lambda)\big) = \mu  M^{{\scriptscriptstyle \mathrm{DES}}, \mu}(\lambda)
\]
and
\[
\lambda \big(1-P_{\mathrm{ab}}^{{\scriptscriptstyle \mathrm{DBS}}, \mu+\theta}(\lambda)\big) = (\mu+\theta)  M^{{\scriptscriptstyle \mathrm{DBS}}, \mu+\theta},(\lambda)
\]
which, using part (a), yield the result.

(c) This follows by reformulating (b) in terms of abandonment rates, noting that
 $\ell^{{\scriptscriptstyle \mathrm{DBS}}, \mu+\theta}(\lambda) = \lambda P_{\mathrm{ab}}^{{\scriptscriptstyle \mathrm{DBS}}, \mu+\theta}(\lambda)$ and
 $\ell^{{\scriptscriptstyle \mathrm{DES}}, \mu}(\lambda) = \lambda P_{\mathrm{ab}}^{{\scriptscriptstyle \mathrm{DES}}, \mu}(\lambda)$.
\end{proof}

Proposition \ref{pro:metricsequiv}(c) reduces evaluation of the abandonment rate in the DES case to the DBS case, addressed in \ref{s:evalelldbs}.
The same applies to evaluation of the derivative of the abandonment rate with respect to the arrival rate $\lambda$, as from Proposition \ref{pro:metricsequiv}(c) we obtain
\begin{equation}
\label{eq:lossrderdes}
(\mu + \theta) \frac{d}{d\lambda} \ell^{{\scriptscriptstyle \mathrm{DES}}, \mu}(\lambda) = 
\theta +
\mu 
\frac{d}{d\lambda} \ell^{{\scriptscriptstyle \mathrm{DBS}}, \mu+\theta}(\lambda).
\end{equation}

\section{Model reformulation as an MARBP and  RB index}
\label{a:rbip}
The MARBP
 concerns the optimal dynamic activation of a collection of stochastic \emph{projects} modeled as RBs ---binary-action (active or passive) MDPs--- subject to given activation constraints.
In the present setting, the ``projects''  are basic nodes, each of which we now view as being fed with its own Poisson job arrival stream with rate $\lambda$ (the total arrival rate) and endowed with an \emph{entry gate} that can be opened (passive action) and shut (active action) by the system controller, as seen in 
Fig.\ \ref{fig:acrrbmodel} (cf.\ Fig.\ \ref{fig:acrmodel}).

\begin{figure}[!htb]
    \centering
    \psfrag{a}{$\lambda$}
    \psfrag{b}[l]{$A_1(t) = 0$}
    \psfrag{c}[l]{$A_n(t) = 1$}
    \psfrag{g}[rb]{$L_1(X_1(t))$}
    \psfrag{h}[lb]{$L_n(X_n(t))$}
    \psfrag{e}[Bb]{$\mu_1$}
    \psfrag{f}[Bb]{$\mu_n$}
    \includegraphics[width=0.65\textwidth]{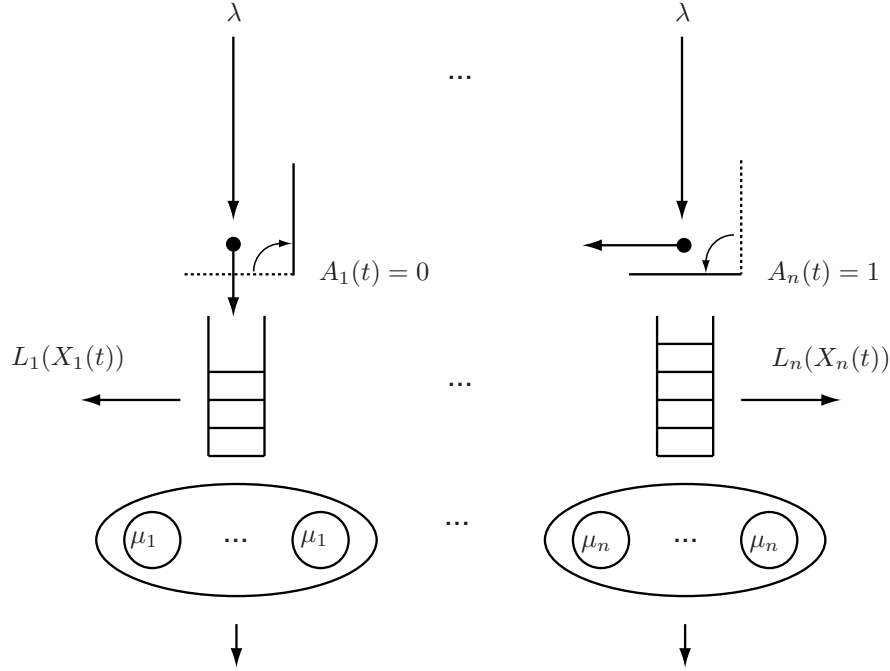}
    \caption{The cloud model viewed as an MARBP.}
    \label{fig:acrrbmodel}
\end{figure}

Instead of modeling actions through the $A(t) \in \{0, 1, \ldots, n\}$ in the formulation in \S \ref{s:md}, 
 we consider binary actions $A_k(t) \in \{0,
1\}$ for basic nodes $k$, where actions $0$ and $1$ mean opening and shutting the node's entry gate, respectively.
Note that $A_k(t) = 0$  means that an arrival at time $t$  would be sent to node $k$, whereas setting $A_k(t) = 1$ for every basic node $k$ (shutting all gates) means that an arrival at time $t$  would be sent to the external node.
 
The requirement that each arrival must be either  outsourced to the external node or routed to a  basic node is thus formulated by 
the constraint that at most one entry gate be open at any time:
 \begin{equation}
 \label{eq:spconstr}
 \sum_{k =1}^n \big(1-A_k(t)\big)  \leqslant 1, \quad t \geqslant 0, 
 \end{equation}
or, equivalently, 
\begin{equation}
\label{eq:2spconstr}
\sum_{k=1}^n A_k(t) \geqslant n-1, \quad t \geqslant 0,
\end{equation}
i.e., the entry gates of at least $n-1$ basic nodes must be shut at any time.
We denote by $\boldsymbol{\Pi}(n-1)$ the class of 
nonanticipative policies satisfying (\ref{eq:2spconstr}).

The total cost rate per unit time incurred 
 when a joint action $\mathbf{a} = (a_k)_{k =1}^n$ satisfying (\ref{eq:2spconstr}) is taken when the system lies in state
$\mathbf{i} = (i_k)_{k =1}^n$  is  
\[
 \sum_{k =1}^n  L_k(i_k) + \lambda C \bigg(1 - \sum_{k =1}^n (1-a_k)\bigg),
\]
which is immediately reformulated to separable form as
\begin{equation}
\label{eq:obj}
-(n-1)  \lambda C + \sum_{k =1}^n \big(L_k(i_k) + \lambda  C a_k\big).
\end{equation}

We can thus reformulate the cloud service provider's cost minimization problem (\ref{eq:panew}) as the MARBP 
\begin{equation}
\label{eq:pa}
\minim_{\boldsymbol{\pi} \in \boldsymbol{\Pi}(n-1)} \, \limsup_{T \to \infty} \frac{1}{T}
\, \Ex_{\mathbf{i}}^{\boldsymbol{\pi}}\bigg[\int_0^T  \sum_{k =1}^n \big(L_k(X_k(t)) + C \lambda  A_k(t)\big) \, dt\bigg],
\end{equation}
where the constant term in (\ref{eq:obj}) is dropped from the objective.

A Lagrangian relaxation approach to (\ref{eq:pa}) leads to consider  
single-node subproblems of the form
\begin{equation}
\label{eq:snpa}
\minim_{\pi_k \in \Pi_k} \, \limsup_{T \to \infty} \frac{1}{T}
\, \Ex_{i_k}^{\pi_k}\bigg[\int_0^T  \big(L_k(X_k(t)) + r \lambda  A_k(t)\big) \, dt\bigg],
\end{equation}
for $k = 1, \ldots, n$, 
where
$\Pi_k$ is the class of nonanticipative \emph{admission control policies} for node $k$  and $r$ is a parameter ranging over $\mathbb{R}$ representing a charge per rejected job.
Thus, (\ref{eq:snpa}) is an average-cost optimal admission control subproblem for node $k$'s M/M/$m_k+$M queue, with the objective of minimizing the average expected rate of deadline-missing and rejection costs.

Let us say that the parametric problem family (\ref{eq:snpa}) is \emph{indexable} if there exists an index $\varphi_k^{\scriptscriptstyle  \mathrm{RB}}\colon \mathbb{Z}_{ +}~\to~\mathbb{R}$   such that, for any value of the rejection charge parameter $r \in \mathbb{R}$, it is optimal in (\ref{eq:snpa}) ---for every initial state $i$--- to reject an arrival finding $X_k(t) = j$ jobs present  (i.e., take the active action $A_k(t) = 1$)  if and only if $\varphi_k^{\scriptscriptstyle  \mathrm{RB}}(j) \geqslant r$, \emph{and} it is optimal to admit it (i.e., take the passive action $A_k(t) = 0$) if and only if $\varphi^{\scriptscriptstyle  \mathrm{RB}}(j) \leqslant r$.
In such a case we call $\varphi_k^{\scriptscriptstyle  \mathrm{RB}}(j)$ node $k$'s \emph{Whittle index}.

Indeed, such a problem family is indexable. 
This follows from \citet[Cor.\ 7.1]{nmmp02}, which, under certain conditions (see \citet[Ass.\ 7.1]{nmmp02}), establishes indexability with a nondecreasing Whittle index in a broader birth-death model.
In the present setting, such conditions reduce to the following: (i) the death rate 
$D_k(i)$ in (\ref{eq:barmu}) is concave nondecreasing in $i_k$; and (ii) the abandonment rate $L_k(i_k)$ is convex nondecreasing in $i_k$.
It is immediate to verify that both conditions hold in the present model.

\citet[pp.\ 396--397]{nmmp02} gives an algorithm for computing the Whittle index for the broader model considered there. In the present setting, the algorithm
reduces to a set of first-order linear recursions  on the index $\varphi^{\scriptscriptstyle  \mathrm{RB}}(i)$ along with auxiliary quantities $z(i+1)$ and $g(i)$, where the node's label $k$ is dropped from the notation. First, let
\begin{align}
\label{eq:phirb0}
\varphi^{\scriptscriptstyle  \mathrm{RB}}(0) & = \frac{L(1)}{D(1)}, 
 & 
z(1) & = 1, & g(0) & = 
  \frac{\lambda D(1)}{\lambda + D(1)}.
\end{align}
Then, for $i = 1, 2, \ldots$, recursively compute
\begin{equation}
\begin{split}
\label{eq:wrcomp}
\varphi^{\scriptscriptstyle  \mathrm{RB}}(i)  & = 
\varphi^{\scriptscriptstyle  \mathrm{RB}}(i-1) + \frac{\displaystyle{\Delta L(i+1) - 
  \varphi^{\scriptscriptstyle  \mathrm{RB}}(i-1)} \Delta D(i+1)}{\displaystyle{\Delta D(i+1) + 
  g(i-1)/\bar{\rho}(i-1)}} \\
  z(i+1) & =  1-\frac{\lambda  D({i})}{\big[\lambda +
D({i})\big]  \big[\lambda + D(i+1)\big] z({i})}, \,
g({i})   = \lambda
\frac{\displaystyle \Delta D(i+1) + 
  \frac{g({i-1})}{\bar{\rho}({i-1})}}{z(i+1) [\lambda + D(i+1)]}.
  \end{split}
\end{equation}
where $D(i)$ is as in (\ref{eq:barmu}), $\bar{\rho}(i) \triangleq
\lambda/D({i+1})$, $\Delta D(i) \triangleq D(i) - D(i-1)$, and similarly with $\Delta L(i)$. 

\begin{remark}
\label{re:rbicases}
\begin{itemize}
\item[\textup{(i)}] Using the above linear recursions  the index values $\varphi^{\scriptscriptstyle  \mathrm{RB}}(0)$, \ldots,
$\varphi^{\scriptscriptstyle  \mathrm{RB}}(i)$ are efficiently computed in $O(i)$ time.
\item[\textup{(ii)}]
Under DBS 
recursion (\ref{eq:wrcomp}) simplifies since $\Delta L(i) = 0$ and $\Delta D(i) = \mu$ for $1 \leqslant i \leqslant m$, whereas
$\Delta L(i) = \Delta D(i) = \theta$ for $i \geqslant m+1$. Hence, 
\begin{equation}
\label{eq:nurbi1m}
\varphi^{\scriptscriptstyle  \mathrm{RB}}(i) = 
\begin{cases}
0,& 0 \leqslant  i \leqslant m-1 \\
\displaystyle \varphi^{\scriptscriptstyle  \mathrm{RB}}(i-1) + \theta \frac{1 - 
  \varphi^{\scriptscriptstyle  \mathrm{RB}}(i-1)}{\displaystyle{\theta + 
  D(i) g(i-1)/\lambda}}, &  i \geqslant m.
\end{cases}
\end{equation}
\item[\textup{(iii)}]
Under DES
recursion  (\ref{eq:wrcomp}) simplifies since $\Delta L(i) = \theta$ and $\Delta D(i) = \mu + \theta$ for $1 \leqslant i \leqslant m$, whereas
$\Delta L(i) = \Delta D(i) = \theta$ for $i \geqslant m+1$. Hence, 
\begin{equation}
\label{eq:nurbi1m2}
\varphi^{\scriptscriptstyle  \mathrm{RB}}(i) = 
\begin{cases}
\theta/(\mu + \theta),& 0 \leqslant  i \leqslant m-1 \\
\displaystyle \varphi^{\scriptscriptstyle  \mathrm{RB}}(i-1) +  \theta \frac{1 - 
  \varphi^{\scriptscriptstyle  \mathrm{RB}}(i-1)}{\displaystyle{\theta + 
  D(i) g(i-1)/\lambda}}, &  i \geqslant m.
\end{cases}
\end{equation}
\end{itemize}
\end{remark}





\end{document}